# REPRESENTATION OF RATIONAL POSITIVE REAL FUNCTIONS OF SEVERAL VARIABLES BY MEANS OF POSITIVE LONG RESOLVENT

M. F. Bessmertnyi

**Abstract.** A rational homogeneous (of degree one) positive real matrix-valued function $f(z)$, $z \in \mathbb{C}^d$ is presented as the Schur complement of a block of the linear pencil $A(z) = z_1 A_1 + \ldots + z_d A_d$ with positive semidefinite matrix coefficients $A_k$. The partial derivative numerators of a rational positive real function are the sums of squares of polynomials.

## 1. Introduction

The long-resolvent Theorem, proved by the author in [3, 4] (see also [5] - [7]), states that every rational $m \times m$ matrix-valued function $f(z) = f(z_1, \cdots, z_d)$ can be represented as the Schur complement

$$f(z) = A_{11}(z) - A_{12}(z) A_{22}(z)^{-1} A_{21}(z), \qquad (1.1)$$

of the block $A_{22}(z)$ of a linear $(m+l) \times (m+l)$ matrix-valued function (linear pencil)

$$A(z) = A_0 + z_1 A_1 + \ldots + z_d A_d = \begin{pmatrix} A_{11}(z) & A_{12}(z) \\ A_{21}(z) & A_{22}(z) \end{pmatrix}. \qquad (1.2)$$

If, moreover, $f(z)$ satisfies the additional condition

(i) $\overline{f(\overline{z}_1, \cdots, \overline{z}_d)} = f(z_1, \cdots, z_d)$,

(ii) $f(z)^T = f(z)$,

(iii) $f(\lambda z_1, \cdots, \lambda z_d) = \lambda f(z_1, \cdots, z_d)$, $\lambda \in \mathbb{C} \setminus \{0\}$),

then one can choose matrices $A_k$, $k = 0, 1, \cdots, d$ to be (i) real (resp. (ii) symmetric, (iii) such that $A_0 = 0$).

Any relation (1.1) is called *a long-resolvent representation*, since for invertible $m \times m$ matrices $f(z)$ representation (1.1) takes the form

$$f(z) = \left[ \pi (A_0 + z_1 A_1 + \ldots + z_d A_d)^{-1} \pi^T \right]^{-1}, \qquad (1.3)$$

where $\pi = (I_m, 0)$ is $m \times (m+l)$ matrix, and $I_m$ is the matrix unit.

There exist many matrix pencil (1.2) representing the same function $f(z)$.

A particular role is a class $\mathbb{RB}_d^{m \times m}$ of $m \times m$ matrix-valued functions (1.1) with *a homogeneous positive real* matrix pencil $A(z) = z_1 A_1 + \ldots + z_d A_d$, $A_k^T = \overline{A}_k = A_k \geq 0$, $k = 1, \cdots, d$, where positive definiteness is understood in the sense of quadratic forms. Functions of class $\mathbb{RB}_d^{m \times m}$, and only they, are characteristic functions of passive electric 2m-poles (see [4], [8] and review [1]), containing ideal transformers and $d$ types of elements (each element of the $k$-th type has an impedance $z_k$).

Suppose $f(z)$ is a function of the class $\mathbb{RB}_d^{m\times m}$. From (1.1) it follows that

$$f(z) = \begin{pmatrix} I_m, & -A_{21}(w)^* A_{22}(w)^{-1*} \end{pmatrix} \begin{pmatrix} A_{11}(z) & A_{12}(z) \\ A_{21}(z) & A_{22}(z) \end{pmatrix} \begin{pmatrix} I_m \\ -A_{22}(z)^{-1} A_{21}(z) \end{pmatrix} = \Phi(w)^* A(z)\, \Phi(z),$$

where $w, z \in \mathbb{C}^d$, $\Phi(w)^* := \overline{\Phi(w)}^T$. Then

$$f(z) + f(z)^* = \sum_{k=1}^d (z_k + \overline{z}_k) \Phi(z)^* A_k\, \Phi(z).$$

Since $A_k \geq 0$, we have

(iv). $f(z) + f(z)^* \geq 0$ for $z \in \Pi^d$,

where $\Pi^d = \{z \in \mathbb{C}^d : \operatorname{Re} z_1 > 0, \cdots, \operatorname{Re} z_d > 0\}$ is open right poly-halfplane.

Function $f(z)$ is a rational. Than from (iv) it follows that

(v). $f(z)$ is holomorphic on $\Pi^d$.

Any matrix-valued function satisfying conditions (i) - (v) is called *a positive real*. The class $\mathbb{RP}_d^{m\times m}$ of *the rational positive real* $m\times m$ matrix-valued functions $f(z_1,\cdots,z_d)$ was introduced by the author in [4] (see also [5]). It's clear that $\mathbb{RB}_d^{m\times m} \subseteq \mathbb{RP}_d^{m\times m}$. The next question is open:

*Are conditions* (i) - (v) *sufficient for belonging rational matrix-valued function $f(z)$ to the class* $\mathbb{RB}_d^{m\times m}$?

Else: *does a matrix-valued function $f(z) \in \mathbb{RP}_d^{m\times m}$ have a long-resolvent representation with positive semidefinite matrices $A_k \geq 0$, $k = 1,\cdots,d$?*

The case of $d = 1$ is not interesting: the classes $\mathbb{RB}_1^{m\times m}$ and $\mathbb{RP}_1^{m\times m}$ coincide and consist of functions of the form $f(z) = z_1 A$ with $m\times m$ matrix $A$ such that $\overline{A} = A^T = A \geq 0$. If $d = 2$, then we also have $\mathbb{RB}_2^{m\times m} = \mathbb{RP}_2^{m\times m}$. If $d \geq 3$, then the question of coincidence of the classes $\mathbb{RB}_d^{m\times m}$ and $\mathbb{RP}_d^{m\times m}$ has remained open till now, with the exception functions of degree 2 and some others.

To solve this question, it is required to obtain a convenient characterization of matrix-valued functions of the class $\mathbb{RB}_d^{m\times m}$. The necessary conditions for belonging $f(z)$ to the class $\mathbb{RB}_d^{m\times m}$ were obtained by the author in [4]: there exist rational real matrix-valued functions $\Phi_k(z)$, $k = 1,\cdots,d$ holomorphic on $\Pi^d$ such that $\Phi_k(\lambda z) = \Phi_k(z)$, $\lambda \in \mathbb{C}\setminus\{0\}$ and

$$f(z) = \sum_{k=1}^d z_k \Phi_k(w)^* \Phi_k(z), \quad w, z \in \mathbb{C}^d. \tag{1.4}$$

The sufficiency of condition (1.4) was proved by Kalyuzhnyi-Verbovetzkyi [17, 18]. He also obtained characterizations of the form (1.4) for various operator generalizations of the class $\mathbb{RB}_d^{m\times m}$, as well as for the operator-valued Schur-Agler functions in the unit polydisk. In [2], for rational matrix-valued Cayley inner Herglotz-Agler functions over the right poly-halfplane was obtained long-resolvent representation (1.1) with matrix pencil (1.2), in which the matrix $A_0$ is skew-Hermitian ($A_0^* = -A_0$) and the other coefficients $A_k$ are positive semidefinite. Thus, the class of Cayley inner rational Herglotz-Agler functions is an extension of the class $\mathbb{R}\mathcal{B}_d^{m \times m}$.

The relation (1.4) requires a representation of non-negative polynomials in the form of a sum of squares of rational functions that are holomorphic in an open right poly-halfplane. Using Artin's solution (see [20]) of Hilbert's 17th problem on the representation of a non-negative rational function of several variables as a sum of squares of rational functions is difficult for two reasons: Artin's theorem says nothing on the location of the singularities of the functions in the decomposition, and, moreover, the proof of the existence of the representation is not constructive. It is interesting to note (J.W. Helton [16]) that each non-commutative "positive" polynomial is the sum of the squares of the polynomials.

In [19] (1968) T. Koga considered a similar class of positive real functions $d$ variables (without condition (iii) homogeneity of degree 1). For a given positive function of several variables, a method for the synthesis of a passive electric circuit was proposed. The T. Koga synthesis is carried out in two stages.

At the first stage, the positive real function $f(z_1, \cdots, z_d)$, which has a degree $n_k$ in each variable $z_k$, is "replaced" by a positive real function $\hat{f}(\varsigma_1, \cdots, \varsigma_n)$ of degree 1 in each variable separately. The variables $\varsigma_1, \cdots, \varsigma_n$, $n = n_1 + \ldots + n_d$ divided into $d$ groups. If the variables of each $k$-th group are replaced by a variable $z_k$, then we get the original positive function $f(z_1, \cdots, z_d)$.

The second stage was based on a modification of Darlington's theorem for functions of several variables. In the univariate case, the Darlington method (see [21]; [14], Chapter V) is based on the representation of a non-negative polynomial as a sum of polynomial squares. T. Koga, as a generalization, suggested

**Koga's Sum-of-Squares Lemma**. *Let $f = f(x_1, \cdots, x_d)$ be a polynomial with real coefficients, quadratic in each variable. If $f \geq 0$ for real $x_i, i = 1, \cdots, d$, then $f = \sum_{j=1}^n h_j^2$, where $h_j$ are polynomials linear in each variable and $n \leq 2^d$.*

As pointed out in 1976 by N.K. Bose [10], Koga's proof is wrong. A counterexample is the non-negative polynomial not representable as a sum of polynomial squares constructed by M.-D. Choi [13]:

$$F(x,y) = x_1^2 y_1^2 + x_2^2 y_2^2 + x_3^2 y_3^2 - 2(x_1 y_1 x_2 y_2 + x_2 y_2 x_3 y_3 + x_1 y_1 x_3 y_3) + 2(x_1^2 y_2^2 + x_2^2 y_3^2 + x_3^2 y_1^2). \quad (1.5)$$

It should be noted that polynomial (1.5) cannot appear when synthesizing a positive real function by the method of T. Koga. In Koga's method the nonnegative polynomial is a partial Wronskian

$$F_k(x) = W_{x_k}[q(x), p(x)] = q(x) \frac{\partial p(x)}{\partial x_k} - p(x) \frac{\partial q(x)}{\partial x_k} \geq 0 \quad (1.6)$$

of a pair of polynomials satisfying the condition $q(z), p(z) \neq 0$ if $z \in \Pi^d$. The polynomial (1.5) cannot be represented in the form (1.6).

The representation (1.6) strongly narrows the class of polynomials under consideration. In some cases, behavior of the rational function $p(z)/q(z)$ allows one to analyze the possibility of representing the Wronskian as a sum of polynomial squares. In particular (see Proposition 8.1), if at least one of the non-negative Wronskians $W_{z_k}[q(z), p(z)]$ cannot be represented as a sum of polynomial squares, then the function $p(z)/q(z)$ cannot be holomorphic on $\Pi^d$.

Hence follows a statement that "rehabilitates" the synthesis of a positive real function according to the method of T. Koga.

**Sum of Squares Theorem**. *If* $P(z)/q(z) \in \mathbb{R}\mathcal{P}_d^{m \times m}$, *then the partial Wronskians*

$$W_{z_k}[q, P] = q(z)\frac{\partial P(z)}{\partial z_k} - P(z)\frac{\partial q(z)}{\partial z_k}, \quad k = 1, \cdots, d \qquad (1.7)$$

*are sums of the squares of polynomials.*

This theorem made it possible to prove the main result of this article: the class $\mathbb{R}\mathcal{P}_d^{m \times m}$ of rational positive real matrix-valued functions coincides with the class $\mathbb{R}\mathcal{B}_d^{m \times m}$ of matrix-valued functions that have a long-resolvent representation with positive semidefinite matrices $\bar{A}_k = A_k^T = A_k \geq 0$, $k = 1, \cdots, d$.

This paper is organized as follows. In Section 2, we explains the terminology and provide preliminary information used in the article.

For the sake of completeness of presentation, in Section 3 we recall the simplest properties of functions of the class $\mathbb{R}\mathcal{P}_d^{m \times m}$ and properties of the degree reduction operator of a rational function. A criterion for the positivity of a real multiaffine function is obtained in Theorem 3.7.

In Section 4, we study Artin's denominator properties of nonnegative forms that cannot be represented as a sum of polynomial squares (PSD not SOS forms). It is proved (Theorem 4.3) that each PSD not SOS form has a minimal Artin's denominator $s(z)$ with irreducible factors do not change sign on $\mathbb{R}^d$.

Theorem 5.1 (Product Polarization Theorem) is a certain "modification" of the long-resolvent Theorem. In particular, Theorem 5.1 provides a convenient representation of the partial Wronskians.

In Section 6, we study the set of the Gram matrices of a given $2n$-form $F$. The Representation Defect Lemma is proved (Lemma 6.8).

A representation of product with one positive semidefinite matrix of the pencil is obtained in Theorem 7.1. The proof is based on Product Polarization Theorem, Representation Defect Lemma and Artin's Theorem.

The sum of squares Theorem is proved in Section 8 (Theorem 8.2). As a preliminary, in Proposition 8.1 we study the singularity points of a function with a nonnegative Wronskian that cannot be represented as a sum of squares of polynomials.

A long-resolvent representation of a rational positive real function with a positive semidefinite matrix pencil is obtained in Theorem 9.1. This proves that the classes of functions $\mathbb{R}\mathcal{P}_d^{m \times m}$ and $\mathbb{R}\mathcal{B}_d^{m \times m}$ coincide.

## 2. Terminology, Designations and Preliminary Information

Let $\mathbb{R}[z]$ be a ring of polynomials in the variables $z_1, \cdots, z_d \in \mathbb{C}$ with real coefficients. We say that $p(z) \in \mathbb{R}[z]$ is *affine* in $z_k$ if $\deg_{z_k} p(z) = 1$, and we say that $p(z)$ is *multiaffine* if it is affine in $z_k$ for all $k = 1, \cdots, d$.

Recall that a *circular region* is a closed or open proper subset of the complex plane, which is bounded by circles (straight lines). In particular, the half-plane is a circular region. We need the following statement about symmetric multiaffine polynomials.

*2.1. Grace-Walsh-Szego Theorem* (see [12], Theorem 2.12). *Let $p$ be a symmetric multiaffine polynomial in $n$ complex variables, let $C$ be an open or closed circular region in $\mathbb{C}$, and let $z_1,\cdots,z_n$ be points in the region $C$. Suppose, further, that either $\deg p = n$ or $C$ is convex (or both). Then there exists at least one point $\xi \in C$ such that $p(z_1,\cdots,z_n) = p(\xi,\cdots,\xi)$.*

Any polynomial $p(z) \in \mathbb{R}[z]$ is called *a form* ($n$-form, a homogeneous polynomial of degree $n$) if $p(\lambda z_1,\cdots,\lambda z_n) = \lambda^n p(z_1,\cdots,z_n)$, $\lambda \in \mathbb{C}$. Everywhere below, with the exception of specially stipulated cases, homogeneous polynomials are considered.

Let $\alpha = (\delta_1,\cdots,\delta_d)$ be a multi-index, where $\delta_j$ are non-negative integers. Then $z^\alpha = z_1^{\delta_1} \cdots z_d^{\delta_d}$ is the monomial of degree $|\alpha| = \delta_1 + ... + \delta_d$.

Let $p_{jk}(z)$ be polynomials. Any matrix $P(z) = \{p_{jk}(z)\}_{j,k=1}^m$ is called *a matrix polynomial*. The polynomial $P(z) = \sum_{j=1}^M B_j z^{\alpha_j}$, where the coefficients $B_j$ are constant matrices, is called *a matrix-valued form* if $\{z^{\alpha_j}\}_{j=1}^M$ is a set of monomials of degree $n$ in variables $z_1,\cdots,z_d$.

Polynomials $P(z) = \{p_{jk}(z)\}_{j,k=1}^m$ and $q(z)$ are said to be *coprime* if for at least one of the polynomials $p_{jk}(z)$ the polynomials $p_{jk}(z)$, $q(z)$ are coprime.

Let $f(z) = \{f_{jk}(z)\}_{j,k=1}^m$ be a rational matrix-valued function. If $q(z)$ is the common denominator of the functions $f_{jk}(z)$, $j,k = 1,\cdots,m$, then $f_{jk}(z) = p_{jk}(z)/q(z)$, where $p_{jk}(z)$ are polynomials. The rational matrix-valued function will be written in the form $f(z) = P(z)/q(z)$, where $P(z) = \{p_{jk}(z)\}_{j,k=1}^m$ is a matrix polynomial and $q(z)$ is a scalar polynomial. In fact, division $P(z)/q(z)$ is the standard operation of multiplying of the matrix $P(z)$ by the number $q(z)^{-1}$.

By definition, put $\partial \{f_{jk}(z)\}_{j,k=1}^m / \partial z_k = \{\partial f_{jk}(z)/\partial z_k\}_{j,k=1}^m$.

In the univariate case, coprime polynomials $p(z)$, $q(z)$ have no common zeros. The situation is different in the case of several variables (the simplest example $z_1/z_2$).

Further, $Z(p) = \{z \in \mathbb{C}^d : p(z) = 0\}$ means the zero set of the polynomial $p(z) \in \mathbb{R}[z_1,\cdots,z_d]$.

*2.2. Theorem* (see [23], Theorem 1.3.2). *Suppose $d > 1$ and $p(z)$, $q(z)$ are coprime polynomials in $d$ complex variables such that $p(0) = q(0) = 0$. If $\Omega$ is a neighborhood of zero in $\mathbb{C}^d$, then:*

(a) *neither of the sets $Z(p) \cap \Omega$ and $Z(q) \cap \Omega$ is a subset of the other;*

(b) *for any $a \in \mathbb{C}$ there exists $z \in \Omega$ such that $q(z) \neq 0$, $p(z)/q(z) = a$.*

In what follows, such a point will be called an *unremovable singularity* of the function $f = p/q$. A rational function $f(z)$ will be called *a regular at a point* $z_0$ if it is defined and takes a finite value at the point $z_0$.

Any matrix $A$ is called *real*, if $\bar{A} = A$, where the bar denotes the replacement of each element of $A$ by a complex conjugate number. The symbol $A^T$ denotes the transposition operation. If $A$ is a matrix with complex elements, then $A^* = \bar{A}^T$ is the Hermitian conjugate matrix.

Any real symmetric $m \times m$ matrix $A$ is called *a positive semidefinite* ($A \geq 0$) if the inequality $\eta^T A \eta \geq 0$ holds for all $\eta \in \mathbb{R}^m$, and *a positive definite* ($A > 0$) if $\eta^T A \eta > 0$ for all $\eta \neq 0$.

A matrix-valued form $F(z)$ will be called *a PSD form* (positive semidefinite) if $F(x) \geq 0$ for $x \in \mathbb{R}^d$. A matrix-valued PSD form will be called *a SOS form* (sum of squares) if $F(z) = H(z) \cdot H(z)^T$, where the elements of the matrix $H(z)$ (possibly rectangular) are scalar forms.

Let $\{z^{\alpha_j}\}_{j=1}^M$ be the set of all monomials of degree $n$ in variables $z_1, \cdots, z_d$. Each $2n$-form $F(z)$ can be represented as

$$F(z) = (z^{\alpha_1}, \cdots, z^{\alpha_M}) \begin{pmatrix} a_{11} & \cdots & a_{1M} \\ \vdots & \ddots & \vdots \\ a_{1M} & \cdots & a_{MM} \end{pmatrix} \begin{pmatrix} z^{\alpha_1} \\ \vdots \\ z^{\alpha_M} \end{pmatrix}. \tag{2.1}$$

The matrix $A = \{a_{jk}\}_{j,k=1}^M$ is called *a Gram matrix* of the $2n$-form $F(z)$. Gram's matrix is not uniquely determined by the $2n$-form. It is known ([22, Theorem 1]) that PSD form $F(z)$ is a SOS form if and only if $F(z)$ has a positive semidefinite Gram matrix.

If $k$ is a field, then $k(x_1, \cdots, x_d)$ denotes the set of rational functions in variables $x_1, \cdots, x_d$ with coefficients from the field $k$.

***2.3. Artin's Theorem*** (see [20, Ch. XI, Corollary 2]). *Let $k$ be a real field admitting only one ordering, and this ordering is Archimedean. Let, further, $f(x) \in k(x_1, \cdots, x_d)$ be a rational function that does not take negative values: $f(a) \geq 0$ for all $a = (a_1, \cdots, a_d) \in k^d$, in which $f(a)$ is defined. Then $f(x)$ there is the sum of squares in $k(x_1, \cdots, x_d)$.* □

Artin's theorem gives a solution to Hilbert's 17th problem. In particular, if $F(z) \in \mathbb{R}[z_1, \cdots, z_d]$ is PSD form, then there exists a form $s(z)$ such that $s(z)^2 F(z)$ is SOS form. The form $s(z)$ will be called *Artin's denominator* of PSD form $F(z)$.

If $F(z)$ is of the SOS form, then $s(z)^2 F(z)$ is also of the SOS form for each form $s(z)$. The question arises: if $F(z)$ is not representable as a sum of squares of forms, then for which forms of $s(z)$ will the form $s(z)^2 F(z)$ also not be SOS form?

***2.4. Proposition*** (see [15, Lemma 2.1]). *Let $F(x)$ be a PSD not SOS form and $s(x)$ an irreducible indefinite form of degree $r$ in $\mathbb{R}[x_1, \cdots, x_d]$. Then $s^2 F$ is also a PSD not SOS form.*

***Proof***. Clearly $s^2 F$ is PSD. If $s^2 F = \sum_k h_k^2$, then for every real tuple $a$ with $s(a) = 0$, it follows that $s^2 F(a) = 0$. This implies $h_k(a)^2 = 0 \ \forall k$ (since $h_k(a)^2$ is PSD), and so on the real variety $s = 0$, we have $h_k = 0$ as well. So (using [9, Theorem 4.5.1]), for each $k$, there exists $g_k$ so that $h_k = s g_k$. This gives $F = \sum_k g_k^2$, which is a contradiction. ∎

**2.5. Corollary**. *Let $F(x)$ be a matrix-valued PSD not SOS form and $s(x)$ an irreducible indefinite form in $\mathbb{R}[x_1,\cdots,x_d]$. Then $s^2F$ is also a PSD not SOS form.*

**Proof**. It suffices to apply the proof of Proposition 2.4 to the diagonal elements of the matrix-valued form $s^2F = H(x)H(x)^T$. ∎

Let $\Pi^d = \{z \in \mathbb{C}^d : \operatorname{Re} z_1 > 0, \cdots, \operatorname{Re} z_d > 0\}$ be an open right poly-halfplane. Any polynomial $p(z) \in \mathbb{R}[z_1,\cdots,z_d]$ is called a *polynomial with the Hurwitz property* [11, 12], or *a stable polynomial* if $p(z) \neq 0$ on $\Pi^d$. Any homogeneous stable polynomial is called *a Hurwitz form*. An example of the Hurwitz form is the nonzero polynomial

$$p(z) = \det(z_1 A_1 + \ldots + z_d A_d), \quad \bar{A}_j = A_j^T = A_j \geq 0, \quad j = 1,\cdots,d. \qquad (2.2)$$

There exist the Hurwitz forms that cannot be represented in the form (2.2) (see [11]).

The polynomial ring is factorial. Than the Hurwitz form is a product $p(z) = \prod [p_k(z)]^{m_k}$ of the irreducible Hurwitz forms. Each irreducible Hurwitz form $p_k(z)$ is indefinite, that is, it changes sign to $\mathbb{R}^d$.

## 3. Positive Real Functions and Degree Reduction Operator

The main object of research in the study of positive real functions is matrix forms of a special type and the possibility of representing such forms as a sum of squares of forms.

**3.1. Proposition.** *If $f(z) = P(z)/q(z) \in \mathbb{R}\mathcal{P}_d^{m\times m}$, then partial Wronskians*

$$W_{z_k}[q(z), P(z)] = q(z)\frac{\partial P(z)}{\partial z_k} - P(z)\frac{\partial q(z)}{\partial z_k}, \quad k = 1,\cdots,d \qquad (3.1)$$

*are PSD forms.*

**Proof**. By the condition $f(z) + f(z)^* \geq 0$ on $\Pi^d$, we put $\lambda = \pm i$. Since $f(\lambda z) = \lambda f(z)$, $\lambda \in \mathbb{C}\setminus\{0\}$, it is observed that

$$\frac{f(z) - f(z)^*}{2i} \geq 0, \operatorname{Im} z_k \geq 0; \quad \frac{f(z) - f(z)^*}{2i} \leq 0, \operatorname{Im} z_k \leq 0, \quad k = 1,\cdots,d. \quad (3.2)$$

Suppose $k = 1$ and $x_2,\cdots,x_d \in \mathbb{R}$; then univariate function $\hat{f}(\varsigma) = f(\varsigma, x_2,\cdots,x_d)$ satisfies the inequality $d\hat{f}(\varsigma)/d\varsigma\big|_{\varsigma\in\mathbb{R}} \geq 0$. From here

$$q(x)\frac{\partial P(x)}{\partial z_1} - P(x)\frac{\partial q(x)}{\partial z_1} = q^2(x)\frac{d\hat{f}}{d\varsigma}(x_1) \geq 0, \quad \text{где } x = (x_1, x_2 \cdots, x_d) \in \mathbb{R}^d. \blacksquare$$

The following statement allows us to eliminate functions $f(z) = P(z)/q(z) \in \mathbb{R}\mathcal{P}_d^{m\times m}$ satisfying the condition

$$\deg_{z_k} P(z) > \deg_{z_k} q(z). \qquad (3.3)$$

**3.2. Proposition.** *Suppose $f(z) = P(z)/q(z) \in \mathbb{R}\mathcal{P}_d^{m\times m}$ satisfies the assumption of (3.3); then there exist a positive semidefinite real $m\times m$ matrix $A_k$ and a form $P_1(z)$ such that:* (a) $\deg_{z_k} P_1(z) = \deg_{z_k} q(z)$; (b) $f_1(z) = P_1(z)/q(z) \in \mathbb{R}\mathcal{P}_d^{m\times m}$; (c) $f(z) = z_k A_k + f_1(z)$.

*Proof.* The degrees of the numerator and denominator of a positive real function $f(z)$ for each variable cannot differ by more than 1. Suppose $k=1$; then there exists a limit

$$\lim_{z_1 \to \infty} \frac{f(z)}{z_1} = A_1(z_2, \cdots, z_d) \neq 0,$$

where the matrix-valued function $A_1(z_2, \cdots, z_d)$ is holomorphic on $\Pi^{d-1}$.

For any $\hat{z} = (z_2, \cdots, z_d) \in \Pi^{d-1}$ the matrix-values function $\hat{f}(z_1) = f(z_1, \hat{z})$ satisfies the inequality $\hat{f}(z_1) + \hat{f}(z_1)^* \geq 0$, $\operatorname{Re} z_1 \geq 0$. The function $\hat{f}(z_1)$ has only simple poles on the imaginary axis (including the pole at $\infty$) and the residues at these poles are positive semidefinite matrices. Then

$$\operatorname*{res}_{z_1 = \infty} \hat{f}(z_1) = \lim_{z_1 \to \infty} \frac{\hat{f}(z_1)}{z_1} = A_1(\hat{z}) \geq 0 \quad \text{for all} \quad \hat{z} \in \Pi^{d-1}.$$

Since $A_1(\hat{z})$ is holomorphic in $\Pi^{d-1}$, we see that $A_1(\hat{z}) \equiv A_1 \geq 0$, where $A_1$ is a constant matrix.

The function $\hat{f}(z_1)$ satisfies the Schwarz-Pick inequality

$$\begin{pmatrix} \dfrac{\hat{f}(z_1) + \hat{f}(z_1)^*}{z_1 + \overline{z}_1} & \dfrac{\hat{f}(z_1) + \hat{f}(x_1)^*}{z_1 + x_1} \\ \dfrac{\hat{f}(z_1)^* + \hat{f}(x_1)}{\overline{z}_1 + x_1} & \dfrac{\hat{f}(x_1) + \hat{f}(x_1)^*}{x_1 + x_1} \end{pmatrix} \geq 0, \tag{3.4}$$

where $\operatorname{Re} z_1 > 0$ and $x_1 > 0$. Limit of (3.4), when $x_1$ is tending to $+\infty$, is a matrix

$$\begin{pmatrix} \dfrac{f(z) + f(z)^*}{z_1 + \overline{z}_1} & A_1 \\ A_1 & A_1 \end{pmatrix} \geq 0. \tag{3.5}$$

Therefore, we have

$$\frac{f(z) + f(z)^*}{z_1 + \overline{z}_1} - A_1 = \begin{pmatrix} I_m & -I_m \end{pmatrix} \begin{pmatrix} \dfrac{f(z) + f(z)^*}{z_1 + \overline{z}_1} & A_1 \\ A_1 & A_1 \end{pmatrix} \begin{pmatrix} I_m \\ -I_m \end{pmatrix} \geq 0. \tag{3.6}$$

It follows from (3.6) that function $f_1(z) = f(z) - A_1 z_1$ satisfies the inequality $f_1(z) + f_1(z)^* \geq 0$, $z \in \Pi^d$. Since $f_1(\lambda z) = \lambda f_1(z)$, $\lambda \in \mathbb{C} \setminus \{0\}$, we see that $f_1(z) \in \mathbb{R}\mathcal{P}_d^{m \times m}$. ∎

Further simplification is based on the use of the degree reduction operator ([19], [12]). In some cases, this allows us to restrict ourselves in considering multi-affine rational positive real functions.

An example of a multi-affine $k$-form is an elementary symmetric polynomial of degree $k$ in variables $\varsigma_1, \cdots, \varsigma_n$:

$$\sigma_k(\varsigma_1, \cdots, \varsigma_n) = \sum_{i_1 < i_2 < \cdots < i_k} \varsigma_{i_1} \varsigma_{i_2} \cdots \varsigma_{i_k}, \quad k = 1, \cdots, n; \quad \sigma_0(\varsigma_1, \cdots, \varsigma_n) = 1. \tag{3.7}$$

Polynomial (3.7) is the sum of $\binom{n}{k} = \dfrac{n!}{k!(n-k)!}$ monomials.

Let $p(z_0, z)$ be an $l$-form, where $(z_0, z) = (z_0, z_1, \cdots, z_d) \in \mathbb{C}^{d+1}$. If $\deg_{z_0} p(z_0, z) = n_0$, then

$$p(z_0, z) = \sum_{k=0}^{n_0} p_{l-k}(z) z_0^k, \qquad (3.8)$$

where $p_{l-k}(z)$ are $(l-k)$-forms in variables $z_1, \cdots, z_d$. Let, further, $\varsigma_1, \cdots, \varsigma_{n_0}$ be new variables.

*3.3. Definition*. Transformation

$$D_{z_0}^{n_0} : \sum_{k=0}^{n_0} p_{l-k}(z) z_0^k \mapsto \sum_{k=0}^{n_0} p_{l-k}(z) \binom{n_0}{k}^{-1} \sigma_k(\varsigma_1, \cdots, \varsigma_{n_0}), \qquad (3.9)$$

is called *a degree reduction operator* in the variable $z_0$. □

If $f(z_0, z) = p(z_0, z) / q(z_0, z)$, where $\deg_{z_0} f(z_0, z) = n_0$, then by definition, put

$$D_{z_0}^{n_0}[f(z_0, z)] := \frac{D_{z_0}^{n_0}[p(z_0, z)]}{D_{z_0}^{n_0}[q(z_0, z)]} = \hat{f}(\varsigma_1, \cdots, \varsigma_{n_0}, z). \square \qquad (3.10)$$

Note that the degrees of $p(z_0, z)$ and $q(z_0, z)$ in the variable $z_0$ may differ, but the same degree reduction operator $D_{z_0}^{n_0}$ is applied to the numerator and denominator.

Since $\binom{n_0}{k}^{-1} \sigma_k(\varsigma_1, \cdots, \varsigma_{n_0}) \Big|_{\varsigma_1 = \cdots = \varsigma_{n_0} = z_0} = z_0^k$, we see that transformations (3.9) and (3.10) are invertible.

It turns out that the degree reduction operator (3.10) preserves the positive reality of the function.

*3.4. Theorem*. If $f(z_0, z) = P(z_0, z) / q(z_0, z) \in \mathbb{R}\mathcal{P}_{d+1}^{m \times m}$, where the polynomials $P(z_0, z)$, $q(z_0, z)$ are coprime and $\deg_{z_0} f(z_0, z) = n_0$, then

$$\hat{f}(\varsigma_1, \cdots, \varsigma_{n_0}, z) = \frac{D_{z_0}^{n_0}[P(z_0, z)]}{D_{z_0}^{n_0}[q(z_0, z)]}$$

is the function of class $\mathbb{R}\mathcal{P}_{d+n_0}^{m \times m}$, affine in each of the variables $\varsigma_1, \cdots, \varsigma_{n_0}$.

We will need a few auxiliary statements.

*3.5. Proposition*. *The coprime numerator and denominator of a scalar positive real function are the Hurwitz forms.*

*Proof*. The homogeneity of the polynomials is obvious. Stability easily follows from a similar fact for functions of one variable having a nonnegative real part in the right half-plane. ∎

*3.6. Lemma*. Let $p(z_0, z)$, $z \in \mathbb{C}^d$ be a Hurwitz form, where $\deg_{z_0} p(z_0, z) = n_0$. Then the polynomial $\hat{p}(\varsigma_1, \cdots, \varsigma_{n_0}, z) = D_{z_0}^{n_0}[p(z_0, z)]$ is also a Hurwitz form.

*Proof*. If $(\hat{z}_1, \cdots, \hat{z}_d) \in \Pi^d$, than a polynomial $p(z_0, \hat{z})$ in one variable $z_0$ has no zeros in the right half-plane. We apply the degree reduction operator in the variable $z_0$. By the Grace-Walsh-Szego Theorem, there exist a point $\xi$, $\operatorname{Re} \xi > 0$ such that $\hat{p}(\varsigma_1, \cdots, \varsigma_n, \hat{z}) = \hat{p}(\xi, \cdots, \xi, \hat{z}) = p(\xi, \hat{z}) \neq 0$. The homogeneity of the polynomial $\hat{p}(\varsigma_1, \cdots, \varsigma_n, \hat{z})$ is obvious. ∎

***Proof of Theorem 3.4***. A matrix-valued function $f(z_0, z)$ is positive real if and only if for any real vector $\xi$ the scalar function $\xi^T f(z_0, z) \xi$ is positive real. Since addition does not deduce from the class of positive functions, we see that

$$z_{d+1} + \frac{p(z_0, z)}{q(z_0, z)} = \frac{z_{d+1} q(z_0, z) + p(z_0, z)}{q(z_0, z)}$$

is a positive real function in variables $z_0, z_1, \cdots, z_d, z_{d+1}$. Its numerator is Hurwitz form. The degree reduction operator is additive. By Lemma 3.6, the polynomial

$$D_{z_0}^{n_0}[z_{d+1} q(z_0, z) + p(z_0, z)] = z_{d+1} D_{z_0}^{n_0}[q(z_0, z)] + D_{z_0}^{n_0}[p(z_0, z)] \tag{3.11}$$

is the Hurwitz form in variables $\varsigma_1, \cdots, \varsigma_{n_0}, z_1, \cdots, z_d, z_{d+1}$. Polynomial (3.11) vanishes at the point $z_{d+1}^0 = -D_{z_0}^{n_0}[p(z_0, z)] / D_{z_0}^{n_0}[q(z_0, z)]$. Since the polynomial (3.11) is stable, we see that $\operatorname{Re} z_{d+1}^0 \leq 0$, $(\varsigma_1, \cdots, \varsigma_{n_0}, z_1, \cdots, z_d) \in \Pi^{d+n_0}$. Hence $\operatorname{Re}\left(D_{z_0}^{n_0}[p(z_0, z)] / D_{z_0}^{n_0}[q(z_0, z)]\right) \geq 0$, $(\varsigma_1, \cdots, \varsigma_{n_0}, z_1, \cdots, z_d) \in \Pi^{d+n_0}$. ∎

***3.7. Theorem*** (**criterion of positivity**). *A rational multiaffine real homogeneous $\left(f(\lambda z) = \lambda f(z),\ \lambda \in \mathbb{C} \setminus \{0\}\right)$ matrix-values function $f(z) = P(z)/q(z)$ belongs to the class $\mathbb{RP}_d^{m \times m}$ iff all Wronskians $W_{z_k}[q, P] = q(z) \frac{\partial P(z)}{\partial z_k} - P(z) \frac{\partial q(z)}{\partial z_k}$, $k = 1, \cdots, d$ are PSD forms.*

***Proof***. The necessity is proved in Statement 3.1. Let us prove the sufficiency. Since $f(z)$ is multiaffine, we see that

$$f(z_k, \hat{z}) = \frac{P(z)}{q(z)} = \frac{z_k P_1(\hat{z}) + P_2(\hat{z})}{z_k q_1(\hat{z}) + q_2(\hat{z})},$$

where $\hat{z} = (z_1, \cdots, z_{k-1}, z_{k+1}, \cdots, z_d)$. If $\hat{z} = \hat{x} \in \mathbb{R}^{d-1}$, than

$$\operatorname{Im} f(z_k, \hat{x}) = \frac{1}{2i} \left( \frac{z_k P_1(\hat{x}) + P_2(\hat{x})}{z_k q_1(\hat{x}) + q_2(\hat{x})} - \frac{\overline{z}_k P_1(\hat{x}) + P_2(\hat{x})}{\overline{z}_k q_1(\hat{x}) + q_2(\hat{x})} \right) =$$

$$= \frac{z_k - \overline{z}_k}{2i} \frac{P_1(\hat{x}) q_2(\hat{x}) - P_2(\hat{x}) q_1(\hat{x})}{|z_k q_1(\hat{x}) + q_2(\hat{x})|^2} = \frac{\operatorname{Im} z_k}{|z_k q_1(\hat{x}) + q_2(\hat{x})|^2} W_{z_k}[q, P](\hat{x}).$$

Hence $\operatorname{Im} f(z_k, \hat{x}) \geq 0$, $\operatorname{Im} z_k > 0$, for each $k = 1, \cdots, d$ (for any real other variables). Therefore (see [5], theorem 2.4) $\operatorname{Im} f(z) \geq 0$ for $z \in i\Pi^d$. Since $f(z)$ is homogeneity, we see that $f(z) \in \mathbb{RP}_d^{m \times m}$. ∎

***3.8. Remark***. *If all Wronskians $W_{z_k}[q, P]$ are PSD forms, and the function $f(z) = P(z)/q(z)$ is not multiaffine, then, generally speaking, it is impossible to guarantee that $f(z)$ belongs to the class $\mathbb{RP}_d^{m \times m}$.*

## 4. Artin's Denominators of PSD not SOS Forms

Let $F(z) \in \mathbb{R}[z_1, \cdots, z_d]$ be a PSD not SOS form. By Artin's theorem, there exists a form $s(z)$ such that $s(z)^2 F(z)$ is a SOS form. The form $s(z)$ is called *Artin's denominator of $F(z)$*.

**4.1. Proposition**. *Suppose $s(z)^2 F(z)$ is a SOS form and each irreducible factor of $s(z)$ is an indefinite form; then $F(z)$ is also a SOS form.*

**Proof**. Assume the converse. Let $F(z)$ be a PSD not SOS form. We decompose the form $s(z)$ into irreducible factors $s(z) = s_1(z) \cdots s_m(z)$, where in the product some factors can occur several times. Consider the forms

$$F_1(z) = s_1^2(z) F(z), \quad F_2(z) = s_2^2(z) F_1(z), \quad \cdots, \quad F_m(z) = s_m^2(z) F_{m-1}(z). \tag{4.1}$$

Successively applying Proposition 2.4 to the forms $F(z), F_1(z), \cdots, F_{m-1}(z)$, we obtain that the form $F_m(z) = s^2(z) F(z)$ is a PSD not SOS form. Contradiction. ∎

**4.2. Definition**. *Artin's denominator $s(z)$ of a PSD not SOS form $F(z)$ is called a minimal Artin denominator if a form $\hat{s}(z) = s(z)/s_j(z)$ is not Artin's denominator of the form $F(z)$ for each of irreducible factor $s_j(z)$ of $s(z)$.* □

**4.3. Theorem**. *Let $F(z)$, $z \in \mathbb{C}^d$ be a PSD not SOS form. Then:*

(a). *There exists a minimal Artin denominator $s(z)$ of $F(z)$;*

(b). *Each irreducible factor of $s(z)$ does not change sign to $\mathbb{R}^d$.*

**Proof**. By Artin's theorem, there exists a form $r(z)$ such that $r(z)^2 F(z)$ there is a SOS form. Each irreducible factor of the form $r(z)$ is either indefinite or does not change sign on $\mathbb{R}^d$. Then $r(z) = s_0(z) s_+(z)$, where all irreducible factors of the form $s_+(z)$ do not change sign on $\mathbb{R}^d$, and the irreducible factors of the form $s_0(z)$ are indefinite. Consider the form $F_1(z) = s_+(z)^2 F(z)$. By condition, $s_0(z)^2 F_1(z) = r(z)^2 F(z)$ is a SOS form. The irreducible factors of the form $s_0(z)$ are indefinite. Then, by Proposition 4.1, $s_+(z)$ is also Artin's denominator of the form $F(z)$. Removing all "excess" irreducible factors from the form $s_+(z)$, we obtain Artin's denominator $s(z)$ with the required properties. ∎

## 5. Product Polarization Theorem

We need a special representation for the product of two forms. Such a representation is a certain "modification" of the long-resolvent Theorem.

**5.1. Theorem** (**Product Polarization Theorem**). *Suppose $q(z)$ is a real form of degree $n$, and $P(z)^T = P(z)$ is a real $m \times m$ matrix-valued form of degree $n+1$, where $z \in \mathbb{C}^d$; then:*

(i). *There exist real symmetric matrices $A_k$, $k = 1, \cdots, d$ such that*

$$q(\varsigma) P(z) = \Psi(\varsigma)(z_1 A_1 + \ldots + z_d A_d) \Psi(z)^T, \quad \varsigma, z \in \mathbb{C}^d, \tag{5.1}$$

*where $\Psi(z) = (z^{\alpha_1} I_m, \cdots, z^{\alpha_N} I_m)$. Here $I_m$ is the matrix unit, and $z^{\alpha_j}$, $j = 1, \cdots, N$ are all monomials of degree $n$ in the variables $z_1, \cdots, z_d$;*

(ii). *The partial Wronskians $W_{z_k}[q(z), P(z)]$ are represented as*

$$q(z) \frac{\partial P(z)}{\partial z_k} - P(z) \frac{\partial q(z)}{\partial z_k} = \Psi(z) A_k \Psi(z)^T, \quad k = 1, \cdots, d. \tag{5.2}$$

The proof of Theorem 5.1 is based on the following statement.

**5.2. Lemma**. *Let* $z^{\alpha_1}, z^{\alpha_2}, \cdots, z^{\alpha_N}$, $N = \dfrac{n!}{d!(n-d)!}$ *be all multiaffine monomials of degree* $n$ *in variables* $z_1, \cdots, z_d$, $1 \le n \le d-1$. *For any multiaffine monomial* $z^\beta$ *of degree* $n+1$ *there exists a linear pencil* $C(z) = z_1 C_1 + \cdots + z_d C_d$ *whose coefficients* $C_j$ $(j = 1, \cdots, d)$ *are real symmetric* $N \times N$ *matrices of rank no more than* $2$ *such that*

$$(z_1 C_1 + \cdots + z_d C_d) \begin{pmatrix} z^{\alpha_1} \\ z^{\alpha_2} \\ \vdots \\ z^{\alpha_N} \end{pmatrix} = \begin{pmatrix} z^\beta \\ 0 \\ \vdots \\ 0 \end{pmatrix}. \tag{5.3}$$

*Proof*. We split a set of variables in $z^\beta$, $z^{\alpha_1}$ into three disjoint subsets:

(i). Common variables of monomials $z^{\alpha_1}$ and $z^\beta$. These variables form the largest common divisor $z^\gamma$ of monomials $z^{\alpha_1}$, $z^\beta$;

(ii). Variables in $z^\beta$ not belonging to $z^{\alpha_1}$. Denote these variables by odd indices: $z_1, z_3, \cdots, z_{2k-1}$. Since $\deg z^\beta = \deg z^{\alpha_1} + 1$, we see that $k \ge 1$;

(iii). Variables in $z^{\alpha_1}$ not belonging to $z^\beta$. Denote these variables by even indices: $z_2, z_4, \cdots, z_{2k-2}$. If $k = 1$, then this subset is empty.

We have

$$z^{\alpha_1} = z_2 z_4 z_6 \cdots z_{2k-4} z_{2k-2} \cdot z^\gamma; \quad z^\beta = z_1 z_3 z_5 \cdots z_{2k-3} z_{2k-1} \cdot z^\gamma.$$

(a). If $k = 1$, then $z^{\alpha_1} = z^\gamma$; $z^\beta = z_1 \cdot z^\gamma = z_1 \cdot z^{\alpha_1}$, then

$$\begin{pmatrix} z_1 & 0 & \cdots & 0 \\ 0 & 0 & \cdots & 0 \\ \vdots & \vdots & \ddots & \vdots \\ 0 & 0 & \cdots & 0 \end{pmatrix} \begin{pmatrix} z^{\alpha_1} \\ z^{\alpha_2} \\ \vdots \\ z^{\alpha_N} \end{pmatrix} = \begin{pmatrix} z^\beta \\ 0 \\ \vdots \\ 0 \end{pmatrix},$$

i.e., $C(z) = z_1 C_1$, where $C_1$ is a real symmetric matrix of rank $1$.

(b). If $k \ge 2$, then we put

$$z^{\alpha_1} = z_2 z_4 z_6 \cdots z_{2k-4} z_{2k-2} \cdot z^\gamma; \qquad z^{\alpha_2} = z_3 z_5 \cdots z_{2k-3} z_{2k-1} \cdot z^\gamma.$$

The monomials $z^{\alpha_j}$, $j = 3, 4, \cdots, 2k-1$ are defined by the formulas

$$z^{\alpha_j} = z^{\alpha_{j-2}} \cdot \dfrac{z_{j-2}}{z_{j-1}}, \quad j = 3, 4, \cdots, 2k-1.$$

If $i \ne j$, then $z^{\alpha_i} \ne z^{\alpha_j}$ and

$$z_1 z^{\alpha_2} = z^\beta, \quad z_j z^{\alpha_j} = z_{j+1} z^{\alpha_{j+2}}, \; j = 1, \cdots, 2k-3, \quad z_{2k-2} z^{\alpha_{2k-2}} = z_{2k-1} z^{\alpha_1}, \quad z_{2k-1} z^{\alpha_{2k-1}} = z^\beta. \tag{5.4}$$

From (5.4) follows the identity

$$\frac{1}{2}\begin{pmatrix} 0 & z_1 & 0 & 0 & 0 & \cdots & 0 & 0 & z_{2k-1} \\ z_1 & 0 & -z_2 & 0 & 0 & \cdots & 0 & 0 & 0 \\ 0 & -z_2 & 0 & z_3 & 0 & \cdots & 0 & 0 & 0 \\ 0 & 0 & z_3 & 0 & -z_4 & \cdots & 0 & 0 & 0 \\ 0 & 0 & 0 & -z_4 & 0 & \cdots & 0 & 0 & 0 \\ \vdots & \vdots & \vdots & \vdots & \vdots & \ddots & \vdots & \vdots & \vdots \\ 0 & 0 & 0 & 0 & 0 & \cdots & 0 & z_{2k-3} & 0 \\ 0 & 0 & 0 & 0 & 0 & \cdots & z_{2k-3} & 0 & -z_{2k-2} \\ z_{2k-1} & 0 & 0 & 0 & 0 & \cdots & 0 & -z_{2k-2} & 0 \end{pmatrix} \begin{pmatrix} z^{\alpha_1} \\ z^{\alpha_2} \\ z^{\alpha_3} \\ z^{\alpha_4} \\ z^{\alpha_5} \\ \vdots \\ z^{\alpha_{2k-3}} \\ z^{\alpha_{2k-2}} \\ z^{\alpha_{2k-1}} \end{pmatrix} = \begin{pmatrix} z^{\beta} \\ 0 \\ 0 \\ 0 \\ 0 \\ \vdots \\ 0 \\ 0 \\ 0 \end{pmatrix}.$$

We have linear pencil

$$C(z) = z_1 C_1 + z_2 C_2 + \cdots + z_{2k-1} C_{2k-1}$$

whose coefficients $C_j$, $j = 1, \cdots, 2k-1$ are a real symmetric matrices of rank $2$. Supplementing the linear pencil $C(z)$ with zeros to the required size, we obtain (5.3). ∎

***Proof of Theorem 5.1***. Let $q(z)$, $P(z)$ be the forms satisfy the conditions

$$\max\{\deg_{z_k} q(z), \deg_{z_k} P(z)\} = n_k, \quad k = 1, \cdots, d.$$

Applying the degree reduction operators $\mathrm{D}_{z_k}^{n_k}$, $k = 1, \cdots, d$ to the forms $q(z)$ and $P(z)$, we obtain the multiaffine forms that satisfy the conditions of the theorem. If Theorem 5.1 is true for the multiaffine forms, then after identifying the corresponding groups of new variables with the variables $z_k$, $k = 1, \cdots, d$, we obtain relation (5.1) for the original forms. Therefore, it suffices to prove Theorem 5.1 for the multiaffine forms.

Let $q(z) = \sum_{j=1}^{N} a_j z^{\alpha_j}$, $P(z) = \sum_{k=1}^{l} B_k z^{\beta_k}$ be the multiaffine forms of degree $n$ and $n+1$. By Lemma 5.2, for a fixed monomial $z^{\alpha_j}$, $j = 1, \cdots, N$, and each monomial $z^{\beta_k}$, $k = 1, \cdots, l$, there exists a symmetric real matrix pencil $C_{jk}(z)$ such that

$$C_{jk}(z) \begin{pmatrix} z^{\alpha_1} \\ \vdots \\ z^{\alpha_{j-1}} \\ z^{\alpha_j} \\ z^{\alpha_{j+1}} \\ \vdots \\ z^{\alpha_N} \end{pmatrix} = \begin{pmatrix} 0 \\ \vdots \\ 0 \\ z^{\beta_k} \\ 0 \\ \vdots \\ 0 \end{pmatrix} - j\text{-th row}, \quad k = 1, \cdots, l, \quad j = 1, \cdots, N.$$

From the coefficients $a_j$, $B_k$ of forms $q(z)$, $P(z)$ and pencils $C_{jk}(z)$ we construct the matrix pencil

$$A(z) = z_1 A_1 + \ldots + z_d A_d = \sum_{j=1}^{N} a_j \left[ \sum_{k=1}^{l} C_{jk}(z) \otimes B_k \right], \tag{5.5}$$

where $C \otimes B = \{c_{jk} \cdot B\}_{j,k=1}^{N}$ is the Kronecker product of the $N \times N$ matrix $C$ and the $m \times m$ matrix $B$. The matrix pencil (5.5) satisfies the condition

$$A(z) \begin{pmatrix} z^{\alpha_1} I_m \\ \vdots \\ z^{\alpha_N} I_m \end{pmatrix} = \begin{pmatrix} a_1 P(z) \\ \vdots \\ a_N P(z) \end{pmatrix}. \qquad (5.6)$$

The properties of reality and symmetry of $A(z)$ are obvious. Multiplying (5.6) on the left by the matrix $(\varsigma^{\alpha_1} I_m, \cdots, \varsigma^{\alpha_N} I_m)$, $\varsigma \in \mathbb{C}^d$, we obtain relation (5.1) for the multiaffine forms. ∎

Note that the matrix pencil $z_1 A_1 + \ldots + z_d A_d$ is not uniquely determined by the forms $q(z)$, $P(z)$. An easy analysis shows that the following statement holds.

**5.3. Proposition**. *Under the conditions of Theorem 5.1, the difference of any two matrix pencils representing the product of forms $q(\varsigma) P(z)$ is a matrix pencil $z_1 S_1 + \ldots + z_d S_d$ such that*

$$(z_1 S_1 + \ldots + z_d S_d) \Psi(z)^T \equiv 0, \quad \Psi(z) S_k \Psi(z)^T \equiv 0. \square \qquad (5.7)$$

## 6. Gram's Matrices and Representation Defect Lemma

Let $F(z)$, $z \in \mathbb{C}^d$ be a real $2n$-form such that $\deg_{z_k} F(z) \leq 2n_k$, $k = 1, \cdots, d$, and let $\mathcal{M} = \{z^{\alpha_j}\}_{j=1}^{N}$ be the set of monomials of degree $n$ satisfying the condition $\deg_{z_k} z^{\alpha_j} \leq n_k \leq n$, $k = 1, \cdots, d$.

Suppose the form $F(z)$ has two Gram's matrices $A_1$, $A_2$; then for the real matrix $S = A_1 - A_2$ the relation holds

$$(z^{\alpha_1}, \cdots, z^{\alpha_N}) \begin{pmatrix} s_{11} & \cdots & s_{1N} \\ \vdots & \ddots & \vdots \\ s_{N1} & \cdots & s_{NN} \end{pmatrix} \begin{pmatrix} z^{\alpha_1} \\ \vdots \\ z^{\alpha_N} \end{pmatrix} = 0, \quad \overline{s}_{ij} = s_{ij} = s_{ji}. \qquad (6.1)$$

The set of matrices $S$ satisfying (6.1) is a linear space $L_0$. Before considering the Representation Defect Lemma, we construct a special basis in the linear space $L_0$.

**6.1. Proposition**. *In the linear space $L_0$ of real symmetric matrices satisfying condition (6.1) there exists a basis such that each basis matrix has either $3$ or $4$ nonzero elements at the intersection of rows and columns corresponding of monomials*

$$z_r^2 z^\gamma, \; z_r z_l z^\gamma, \; z_l^2 z^\gamma, \qquad (6.2)$$

$$z_r z^{\gamma_1}; \; z_l z^{\gamma_1}; \; z_l z^{\gamma_2}; \; z_r z^{\gamma_2}, \; z^{\gamma_1} \neq z^{\gamma_2}. \qquad (6.3)$$

*The nonzero submatrices of the corresponding basis matrices are determined by the relations*

$$(z_r^2 z^\gamma, \; z_r z_l z^\gamma, \; z_l^2 z^\gamma) \begin{pmatrix} 0 & 0 & -1 \\ 0 & 2 & 0 \\ -1 & 0 & 0 \end{pmatrix} \begin{pmatrix} z_r^2 z^\gamma \\ z_r z_l z^\gamma \\ z_l^2 z^\gamma \end{pmatrix} \equiv 0, \qquad (6.4)$$

$$(z_r z^{\gamma_1},\ z_l z^{\gamma_1},\ z_l z^{\gamma_2},\ z_r z^{\gamma_2}) \begin{pmatrix} 0 & 0 & 1 & 0 \\ 0 & 0 & 0 & -1 \\ 1 & 0 & 0 & 0 \\ 0 & -1 & 0 & 0 \end{pmatrix} \begin{pmatrix} z_r z^{\gamma_1} \\ z_l z^{\gamma_1} \\ z_l z^{\gamma_2} \\ z_r z^{\gamma_2} \end{pmatrix} \equiv 0. \quad (6.5)$$

**6.2. Comment**. Relation (6.1) can be conveniently rewritten as

$$F(z) = \sum_{i,j=1}^{N} s_{ij} z^{\alpha_i} \cdot z^{\alpha_j} = \sum_k c_k z^{\beta_k} = 0, \quad (6.6)$$

where $\deg_{z_k} z^{\beta_j} \leq 2n_k$, $k = 1, \cdots, d$, and $z^{\beta_i} \neq z^{\beta_j}$, $i \neq j$.

It follows from (6.6) that the elements $s_{ij}$ of the matrix (6.1) satisfy the conditions

$$\sum_{\alpha_i + \alpha_j = \beta_k} s_{ij} = 0. \square \quad (6.7)$$

We will need several lemmas.

Let $\Pi_\beta$ be the set of all unordered pairs $\pi = (z^{\alpha_i}, z^{\alpha_j})$, $z^{\alpha_i}, z^{\alpha_j} \in \mathcal{M}$ such that $z^{\alpha_i} z^{\alpha_j} = z^\beta$, where $z^\beta = z_1^{r_1} \cdots z_d^{r_d}$ is a fixed monomial of degree $2n$.

**6.3. Lemma**. Let $\beta = (r_1, \cdots, r_d)$ be a multi-index. Monomial $z^\alpha = z_1^{\delta_1} \cdots z_d^{\delta_d} \in \mathcal{M}$ belongs to some pair $\pi \in \Pi_\beta$ iff

$$\max\{(r_k - n_k), 0\} \leq \delta_k \leq \min\{r_k, n_k\}, \quad k = 1, \cdots, d. \quad (6.8)$$

**Proof**. If $z^\alpha z^{\alpha_1} = z^\beta$, where $z^\alpha = z_1^{\delta_1} \cdots z_d^{\delta_d}$, $z^{\alpha_1} = z_1^{v_1} \cdots z_d^{v_d} \in \mathcal{M}$, then $v_k = r_k - \delta_k$. Hence

$$0 \leq \delta_k \leq n_k, \qquad r_k - n_k \leq \delta_k \leq r_k.$$

These inequalities coincide with (6.8). Sufficiency is obvious. ∎

**6.4. Definition**. If a monomial $z^\alpha = z_1^{\delta_1} \cdots z_r^{\delta_r} \cdots z_l^{\delta_l} \cdots z_d^{\delta_d} \in \mathcal{M}$ satisfies conditions $\delta_r < n_r$, $\delta_l \geq 1$, then a transformation

$$z_1^{\delta_1} \cdots z_r^{\delta_r} \cdots z_l^{\delta_l} \cdots z_d^{\delta_d} \mapsto z_1^{\delta_1} \cdots z_r^{\delta_r + 1} \cdots z_l^{\delta_l - 1} \cdots z_d^{\delta_d} \quad (6.9)$$

is called *an elementary monomial transformation*. $\square$

**6.5. Lemma**. Suppose monomials $z^{\alpha_i}, z^{\alpha_j}$, $z^{\alpha_i} \neq z^{\alpha_j}$ satisfies the assumption (6.8); then:

(i). *There exists a chain of pairwise different monomials satisfying conditions* (6.8) *"connecting" monomials $z^{\alpha_i}$ and $z^{\alpha_j}$;*

(ii). *Any two neighboring monomials of the chain are related by an elementary transformation*.

**Proof**. Let $z^{\alpha_i} = z_1^{\delta_1} \cdots z_d^{\delta_d}$ and $z^{\alpha_j} = z_1^{v_1} \cdots z_d^{v_d}$ be two different monomials satisfying condition (6.8). If $k_j = v_j - \delta_j$, $j = 1, \cdots, d$, then

$$-\min\{r_j, n_j\} \leq k_j \leq \min\{r_j, n_j\} \quad \text{and} \quad \sum_{j=1}^{d} k_j = 0.$$

We represent an integer tuple $(k_1,\cdots,k_d)$ in the form of a componentwise sum of the elementary tuples. Each elementary tuple $(m_1,\cdots,m_d)$ contains only two nonzero components $+1$ and $-1$. The numbers $+1$ are associated with $k_j > 0$, and $-1$ with $k_j < 0$. The sequence of the elementary tuples $(m_1,\cdots,m_d)$ corresponds to the sequence of elementary transformations of monomials. We choose as the initial monomial the monomial $z^{\alpha_i} = z_1^{\delta_1} \cdots z_d^{\delta_d}$. At each step we obtain a monomial that satisfies conditions (6.8) and does not coincide with the monomials obtained in the previous steps. At the last step, we get the monomial $z^{\alpha_j} = z_1^{\nu_1} \cdots z_d^{\nu_d}$. ∎

We say that elements $\pi_1, \pi_2 \in \Pi_\beta$ are *connected by an elementary transformation* if one of the monomials of the pair $\pi_1$ is connected by an elementary transformation with one of the monomials of the pair $\pi_2$.

Note that any element of $(z^{\alpha_k}, z^{\alpha_v}) \in \Pi_\beta$ is uniquely determined by the choice of one of the monomials of the pair.

**6.6. Corollary.** *If $\pi_1, \pi_2 \in \Pi_\beta$ and $\pi_1 \neq \pi_2$, then:*

(i). *The set $\Pi_\beta$ contains a chain of elements that "connects" $\pi_1$ and $\pi_2$;*

(ii). *Any two neighboring elements of the chain are related by an elementary transformation.*

**Proof.** It suffices to consider the chain connecting the first monomials from $\pi_1$, $\pi_2$. ∎

**6.7. Lemma.** *If $\Pi_\beta$ contain $m \geq 2$ elements, then in the set $\Pi_\beta$ there exist $(m-1)$ different pairs $\{\pi_i, \pi_j\}$ such that the elements of each pair are connected by an elementary transformation.*

**Proof.** Let us assign a finite graph to the set $\Pi_\beta$. The vertices of the graph are elements of the set $\Pi_\beta$. The edges of the graph form pairs $\{\pi_i, \pi_j\}$ of elements connected by an elementary transformation. By Corollary 6.6, the graph is connected. The graph tree (that is, a connected subgraph with the same vertices, but without cycles) contains the number of edges one less than the number of vertices. Since $m \geq 2$ and in the graph tree different edges are incident to different pairs of vertices, the number of different pairs $\{\pi_i, \pi_j\}$ connected by an elementary transformation is $m-1$. ∎

**Proof of Proposition 6.1.** If each $\Pi_{\beta_k}$ contains at most one element, then $\dim L_0 = 0$. Suppose some of the sets $\Pi_{\beta_k}$ contain $m \geq 2$ elements; then $m-1$ coefficients $s_{ij}$, $i \leq j$ in the sum (6.7) can be chosen arbitrary, and the remaining coefficient is determined from the condition $c_k = 0$. Thus, such set $\Pi_{\beta_k}$ defines a subspace of dimension $m-1$ of the space $L_0$. Let us show that in each such subspace there exists a basis formed by the matrices of the form (6.4) and (6.5).

By Lemma 6.7, in $\Pi_{\beta_k}$ there exist $m-1$ different pairs $\{\pi_i, \pi_j\}$ of elements such that $\pi_i$ and $\pi_j$ are connected by an elementary transformation. The following cases are possible:

(a). One of the elements of a pair $\{\pi_1, \pi_2\}$, say $\pi_1$, has the form $\pi_1 = (z^{\alpha_i}, z^{\alpha_i})$, that is $(z^{\alpha_i})^2 = z^{\beta_k}$. Let the second element of the pair be $\pi_2 = (z^{\alpha_j}, z^{\alpha_l})$ and the monomials $z^{\alpha_i}$, $z^{\alpha_j}$ are connected by an elementary transformation. Then there exist a variables $z_r, z_l$ such that

$$z^{\alpha_j} = z_1^{\delta_1} \cdots z_r^{\delta_r+1} \cdots z_l^{\delta_l-1} \cdots z_d^{\delta_d}, \quad z^{\alpha_i} = z_1^{\delta_1} \cdots z_r^{\delta_r} \cdots z_l^{\delta_l} \cdots z_d^{\delta_d}, \quad z^{\alpha_l} = z_1^{\delta_1} \cdots z_r^{\delta_r-1} \cdots z_l^{\delta_l+1} \cdots z_d^{\delta_d},$$

or

$$z^{\alpha_j} = z_r^2 z^\gamma, \quad z^{\alpha_i} = z_r z_l z^\gamma, \quad z^{\alpha_l} = z_l^2 z^\gamma, \quad z^\gamma = z_1^{\delta_1} \cdots z_r^{\delta_r-1} \cdots z_l^{\delta_l-1} \cdots z_d^{\delta_d}.$$

This triplet of monomials corresponds to the basis matrix (6.4).

(b). $\pi_1 = (z^{\alpha_i}, z^{\alpha_j})$, $z^{\alpha_i} \neq z^{\alpha_j}$; $\pi_2 = (z^{\alpha_l}, z^{\alpha_s})$, $z^{\alpha_l} \neq z^{\alpha_s}$, and monomials $z^{\alpha_i}$, $z^{\alpha_l}$ are connected by an elementary transformation. Then there exist variables $z_r, z_l$ such that

$$z^{\alpha_i} = z_1^{\delta_1} \cdots z_r^{\delta_r} \cdots z_l^{\delta_l} \cdots z_d^{\delta_d}, \quad z^{\alpha_l} = z_1^{\delta_1} \cdots z_r^{\delta_r-1} \cdots z_l^{\delta_l+1} \cdots z_d^{\delta_d},$$

$$z^{\alpha_j} = z_1^{\beta_1-\delta_1} \cdots z_r^{\beta_r-\delta_r} \cdots z_l^{\beta_l-\delta_l} \cdots z_d^{\beta_d-\delta_d}, \quad z^{\alpha_s} = z_1^{\beta_1-\delta_1} \cdots z_r^{\beta_r-\delta_r+1} \cdots z_l^{\beta_l-\delta_l-1} \cdots z_d^{\beta_d-\delta_d},$$

or

$$z^{\alpha_i} = z_r z^{\gamma_1}, \quad z^{\alpha_l} = z_l z^{\gamma_1}, \quad z^{\alpha_j} = z_l z^{\gamma_2}, \quad z^{\alpha_s} = z_r z^{\gamma_2}, \tag{6.9}$$

where $z^{\gamma_1} = z_1^{\delta_1} \cdots z_r^{\delta_r-1} \cdots z_l^{\delta_l} \cdots z_d^{\delta_d}$, $z^{\gamma_2} = z_1^{\eta_1-\delta_1} \cdots z_r^{r_r-\delta_r} \cdots z_l^{\eta-\delta_l-1} \cdots z_d^{r_d-\delta_d}$.

In addition $z^{\gamma_1} \neq z^{\gamma_2}$. Otherwise, if $z^{\gamma_1} = z^{\gamma_2}$, then $\pi_1 = \pi_2$. Contradiction. The monomials (6.9) define the basis matrix (6.5).

Since all pairs $\{\pi_i, \pi_j\}$ are different, then corresponding $(m-1)$ matrices are linearly independent. ∎

**6.8. Representation Defect Lemma.** *Let $z^{\alpha_j}$, $j=1,\cdots,N$ be all monomials of degree $n$ in variables $z_1,\cdots,z_d$ such that $\deg_{z_k} z^{\alpha_j} \leq n_k$, $k=1,\cdots,d$.*

*If a symmetric real matrix $S_d$ satisfy the conditions*

$$(z^{\alpha_1},\cdots z^{\alpha_N})S_d \begin{pmatrix} z^{\alpha_1} \\ \vdots \\ z^{\alpha_N} \end{pmatrix} \equiv 0, \quad S_d \begin{pmatrix} \partial^{n_d} z^{\alpha_1}/\partial z_d^{n_d} \\ \vdots \\ \partial^{n_d} z^{\alpha_N}/\partial z_d^{n_d} \end{pmatrix} \equiv 0, \tag{6.10}$$

*then there exist the symmetric real matrices $S_k$, $k=1,\cdots,(d-1)$ for which*

$$(z_1 S_1 + \ldots + z_{d-1} S_{d-1} + z_d S_d)\begin{pmatrix} z^{\alpha_1} \\ \vdots \\ z^{\alpha_N} \end{pmatrix} = \begin{pmatrix} 0 \\ \vdots \\ 0 \end{pmatrix}. \tag{6.11}$$

*Proof.* Without loss of generality, we can assume that

$$(z^{\alpha_1},\cdots,z^{\alpha_N}) = (z_d^{n_d} \varphi^T(\hat{z}), \; \psi^T(\hat{z},z_d)),$$

where $\deg_{z_d} \psi(\hat{z},z_d) \leq n_d - 1$, and the components of $\varphi(\hat{z})$ are different monomials in the variables $(z_1,\cdots,z_{d-1}) = \hat{z}$. Dividing the matrix $S_d$ into blocks, from the second relation (6.10) we obtain

$$\begin{pmatrix} S_{11} & S_{12} \\ S_{12}^T & \hat{S}_d \end{pmatrix} \begin{pmatrix} n_d! \varphi(\hat{z}) \\ 0 \end{pmatrix} \equiv 0.$$

Hence $S_{11} = 0$, $S_{12}^T = 0$. Then $S_d = \begin{pmatrix} 0 & 0 \\ 0 & \hat{S}_d \end{pmatrix}$ and $\psi^T(\hat{z}, z_d)\hat{S}_d\,\psi(\hat{z}, z_d)) \equiv 0$.

We rewrite (6.11) in block form

$$\begin{pmatrix} A_{11}(\hat{z}) & A_{12}(\hat{z}) \\ A_{12}(\hat{z})^T & A_{22}(\hat{z}) + z_d\hat{S}_d \end{pmatrix}\begin{pmatrix} z_d^{n_d}\varphi(z) \\ \psi(\hat{z}, z_d) \end{pmatrix} = \begin{pmatrix} 0 \\ 0 \end{pmatrix}. \tag{6.12}$$

Matrix $\hat{S}_d$ is a linear combination of basis matrices $\hat{S}_{d,j}$ of the form (6.4), (6.5). If there exists a solution $A_{11}(\hat{z}), A_{12}(\hat{z}), A_{22}(\hat{z})$ of equation (6.12) for each basis matrix $\hat{S}_{d,j}$, then the general solution (6.12) will be a linear combination of basic solutions.

The basis matrix (6.4) is determined by the monomials $z_r^2 z^\gamma$, $z_r z_l z^\gamma$, $z_l^2 z^\gamma$. By construction, we have $\deg_{z_d} \psi(\hat{z}, z_d) \leq n_d - 1$. Among the components of the vector $(z_d^{n_d}\varphi^T(\hat{z}), \psi^T(\hat{z}, z_d))$ there exist monomials $z_d z_r z^\gamma$ and $z_d z_l z^\gamma$. The solution of the equation (6.12) for the basis matrix (6.4) is the matrix

$$\begin{pmatrix} 0 & 0 & 0 & -z_l & z_r \\ 0 & 0 & z_l & -z_r & 0 \\ 0 & z_l & 0 & 0 & -z_d \\ -z_l & -z_r & 0 & 2z_d & 0 \\ z_r & 0 & -z_d & 0 & 0 \end{pmatrix}\begin{pmatrix} z_d z_r z^\gamma \\ z_d z_l z^\gamma \\ z_r^2 z^\gamma \\ z_r z_l z^\gamma \\ z_l^2 z^\gamma \end{pmatrix} \equiv 0.$$

For the monomials $z_r z^{\gamma_1}$, $z_l z^{\gamma_1}$, $z_l z^{\gamma_2}$, $z_r z^{\gamma_2}$ defining the basis matrix (6.5), there exist monomials $z_d z^{\gamma_1}$ and $z_d z^{\gamma_2}$ such that

$$\begin{pmatrix} 0 & 0 & -z_v & z_k & 0 & 0 \\ 0 & 0 & 0 & 0 & -z_k & z_v \\ -z_v & 0 & 0 & 0 & z_d & 0 \\ z_k & 0 & 0 & 0 & 0 & -z_d \\ 0 & -z_k & z_d & 0 & 0 & 0 \\ 0 & z_v & 0 & -z_d & 0 & 0 \end{pmatrix}\begin{pmatrix} z_d z^{\gamma_2} \\ z_d z^{\gamma_1} \\ z_k z^{\gamma_1} \\ z_v z^{\gamma_1} \\ z_v z^{\gamma_2} \\ z_k z^{\gamma_2} \end{pmatrix} \equiv 0.$$

The general solution of equation (6.12) is a linear combination of such solutions. ∎

## 7. Rational Functions with PSD not SOS Wronskian

It follows from Theorem 5.1 that rational function $f(z) = p(z)/q(z)$ is represented as

$$f(z) = \frac{p(z)}{q(z)} = \frac{\Psi(\varsigma)}{q(\varsigma)}(z_1 A_1 + z_2 A_2 + \ldots + z_d A_d)\frac{\Psi(z)^T}{q(z)}, \quad \varsigma, z \in \mathbb{C}^d. \tag{7.1}$$

Let $f(z) = p(z)/q(z)$ be a rational function with PSD partial Wronskian $W_{z_k}[q, p]$. Let us investigate the possibility of obtaining a representation of the form (7.1) with a positive semidefinite matrix $A_k \geq 0$.

**7.1. Theorem**. Let $q(z)$, $p(z)$, $z \in \mathbb{C}^d$ be real forms, where $\deg p(z) = 1 + \deg q(z)$, $\deg_{z_1} p(z) = \deg_{z_1} q(z)$. If partial Wronskian $W_{z_1}[q, p]$ is a PSD not SOS form with Artin's denominator $s(z)$, then:

(i). There exists a matrix pencil $A(z) = z_1 A_1 + z_2 A_2 + ... + z_d A_d$ with a positive semidefinite matrix $A_1 \geq 0$ such that

$$f(z) = \frac{p(z)}{q(z)} = \frac{\Psi(\varsigma)}{q(\varsigma)s(\varsigma)}(z_1 A_1 + z_2 A_2 + ... + z_d A_d)\frac{\Psi(z)^T}{q(z)s(z)}, \qquad \varsigma, z \in \mathbb{C}^d, \tag{7.2}$$

where $\Psi(z)$ is a row vector of monomials $z^{\alpha_j}$ in the variables $z_1, \cdots, z_d$ satisfying the conditions $\deg z^{\alpha_j} = \deg q(z)s(z)$, $\deg_{z_1} \Psi(z) = \deg_{z_1} q(z)s(z) = n_1$;

(ii). Partial Wronskians are represented as

$$W_k[q(z), p(z)] = \frac{\Psi(z)}{s(z)} A_k \frac{\Psi(z)^T}{s(z)}, \qquad k = 1, \cdots, d. \tag{7.3}$$

**Proof**. By Theorem 5.1, there exists a real symmetric matrix pencil $z_1 B_1 + z_2 B_2 + ... + z_d B_d$ such that

$$q(\varsigma)s(\varsigma)p(z)s(z) = \Psi(\varsigma)[z_1 B_1 + z_2 B_2 + ... + z_d B_d]\Psi(z)^T, \quad \varsigma, z \in \mathbb{C}^d, \tag{7.4}$$

$$W_{z_k}[qs, ps] = \Psi(z) B_k \Psi(z)^T, \quad k = 1, \cdots, d. \tag{7.5}$$

Differentiating (7.4) $n_1 + 1$ times with respect to the variable $z_1$, we obtain

$$q(\varsigma)s(\varsigma)\frac{\partial^{n_1+1} p(z)s(z)}{\partial z_1^{n_1+1}} = (n_1 + 1)\Psi(\varsigma) B_1 \frac{\partial^{n_1}\Psi(z)^T}{\partial z_1^{n_1}} + \Psi(\varsigma)(z_1 B_1 + ... + z_d B_d)\frac{\partial^{n_1+1}\Psi(z)^T}{\partial z_1^{n_1+1}}. \tag{7.6}$$

By condition, we have $\dfrac{\partial^{n_1+1} p(z)s(z)}{\partial z_1^{n_1+1}} \equiv 0$, $\dfrac{\partial^{n_1+1}\Psi(z)^T}{\partial z_1^{n_1+1}} \equiv 0$. Using (7.6), we obtain

$$B_1 \frac{\partial^{n_1}\Psi(z)^T}{\partial z_1^{n_1}} \equiv 0. \tag{7.7}$$

The form

$$\Psi(z) B_1 \Psi(z)^T = W_{z_1}[qs, ps] = qs\frac{\partial(ps)}{\partial z_1} - ps\frac{\partial(qs)}{\partial z_1} = s^2\left(q\frac{\partial p}{\partial z_1} - p\frac{\partial q}{\partial z_1}\right) = s^2 W_{z_1}[q, p]$$

is SOS form. The matrix $B_1$ is its Gram matrix. Among the Gram matrices of the SOS form, there exists a positive semidefinite matrix $A_1 \geq 0$:

$$W_{z_1}[qs, ps] = \Psi(z) B_1 \Psi(z)^T = \Psi(z) A_1 \Psi(z)^T.$$

The degree of the SOS form $W_{z_1}[qs, ps]$ (and hence PSD) with respect to the variable $z_1$ is necessarily even. Then

$$\deg_{z_1} \Psi(z) A_1 \Psi(z)^T = \deg_{z_1} W_{z_1}[qs, ps] \leq 2n_1 - 2. \tag{7.8}$$

It follows from (7.8) that diagonal elements $a_{jj}$ of the matrix $A_1$ corresponding to the monomials $z^{\alpha_j}$ of the vector $\Psi(z)^T$ containing $z_1^{n_1}$ are necessarily equal to $0$. From the positive semidefiniteness of $A_1$ we obtain

$$A_1 \frac{\partial^{n_1} \Psi(z)^T}{\partial z_1^{n_1}} \equiv 0. \tag{7.9}$$

By (7.9) and (7.7) it follows that matrix $S_1 = A_1 - B_1$, $S_1^T = S_1 = \overline{S}_1$ satisfies the conditions

$$\Psi(z) S_1 \Psi(z)^T \equiv 0, \qquad S_1 \frac{\partial^{n_1} \Psi(z)^T}{\partial z_1^{n_1}} \equiv 0.$$

Using Representation Defect Lemma, we get

$$\Psi(\varsigma)(z_1 S_1 + z_2 S_2 + ... + z_d S_d) \Psi(z)^T \equiv 0, \qquad \varsigma, z \in \mathbb{C}^d, \tag{7.10}$$

where $S_k$, $k = 2, \cdots, d$ are the real symmetric matrices. Adding (7.10) to (7.4) and dividing both sides of the resulting identity by the product $q(\varsigma) s(\varsigma) q(z) s(z)$, we obtain (7.2).

The relations (7.3) follow from the identities $s^2 W_{z_k}[q, p] = W_{z_k}[qs, ps] = \Psi(z) A_k \Psi(z)^T$. ∎

## 8. Sum of Squares Theorem

By $\mathcal{H}^d = i\Pi^d = \{z \in \mathbb{C}^d : \operatorname{Im} z_k > 0, \ k = 1, \cdots, d\}$ denote the open upper poli-halfplane. The closure of $\mathcal{H}^d$ is $\overline{\mathcal{H}}^d = \{z \in \mathbb{C}^d : \operatorname{Im} z_k \geq 0, \ k = 1, \cdots, d\}$. To prove the Sum of Squares Theorem, we need the following statement.

**8.1. Proposition**. *Let $p(z), q(z) \in \mathbb{R}[z_1, \cdots, z_d]$, $d \geq 3$ be are forms such that:*

(i). $\deg p(z) = 1 + \deg q(z)$;

(ii). $\deg_{z_1} p(z) = \deg_{z_1} q(z)$;

(iii). *The partial Wronskian $W_{z_1}[q, p] = q(z) \dfrac{\partial p(z)}{\partial z_1} - p(z) \dfrac{\partial q(z)}{\partial z_1}$ is a PSD not SOS form.*

*Then $f(z) = p(z)/q(z)$ has a singularity at some point $(z_1', z_2' \cdots, z_d') \in \overline{\mathcal{H}}^d$, $\operatorname{Im} z_1' > 0$.*

First, we prove

**8.2. Лемма**. *Let $s_0(z) \in \mathbb{R}[z_1, \cdots, z_d]$, $\partial s_0(z)/\partial z_1 \not\equiv 0$ be a irreducible form that does not change sign to $\mathbb{R}^d$. If the form $s_0(z)$ is not a divisor of the polynomial $h(z) \in \mathbb{R}[z_1, \cdots, z_d]$, then there exists a point $z' = (\eta_1', x_2', \cdots, x_d')$, $\operatorname{Im} \eta_1' > 0$, $x_2', \cdots, x_d' \in \mathbb{R}$ such that $s_0(z') = 0$, $h(z') \neq 0$.*

***Proof***. Each irreducible form that does not change sign on $\mathbb{R}^d$ depends on at least two variables. Since $\partial s_0(z)/\partial z_1 \not\equiv 0$, we see that $\deg_{z_1} s_0(z) \geq 2$ and there exists a variable $z_k$, $k \neq 1$ for which $\partial s_0(z)/\partial z_k \not\equiv 0$.

The form $s_0(z)$ is not a divisor of the polynomial $h(z)$. Then there exist $x_2', \cdots, x_d' \in \mathbb{R}$, $\delta > 0$ such that for any $x_2, \cdots, x_d \in \mathbb{R}$ satisfying the condition $|x_k - x_k'| < \delta$, $k = 2, \cdots, d$, at least one zero of the polynomial $s_0(z_1, x_2, \cdots, x_d)$ in the variable $z_1$ is not a zero of the polynomial $h(z_1, x_2, \cdots, x_d)$. The real zeros of $s_0(z_1, x_2, \cdots, x_d)$ have even multiplicity, and the remaining zeros form complex conjugate pairs. For an arbitrarily small change in the numbers $x_2, \cdots, x_d \in \mathbb{R}$ (and hence the coefficients of the polynomial), all real multiple zeros split into complex conjugate pairs. Then there exists $\eta_1'$, $\operatorname{Im}\eta_1' > 0$ such that $s_0(\eta_1', x_2', \cdots, x_d') = 0$ and $h(\eta_1', x_2', \cdots, x_d') \neq 0$. ∎

***Proof of Proposition 8.1***. By Theorem 4.3, for PSD not SOS form $W_{z_1}[q, p]$, there exists a minimal Artin's denominator $s(z) = s_1(z)^{k_1} \cdots s_n(z)^{k_n}$ with irreducible factors $s_j(z)$, $j = 1, \cdots, n$ that do not change sign to $\mathbb{R}^d$. There are two alternative possibilities for $s(z)$:

(a). There exists an irreducible factor $s_j(z)$ depending on the variable $z_1$;

(b). Artin's denominator $s(z)$ does not depend on the variable $z_1$.

By Theorem 7.1, there exists a symmetric matrix pencil $A(z) = z_1 A_1 + z_2 A_2 + \ldots + z_d A_d$ with a positive semidefinite matrix $A_1 \geq 0$ such that

$$f(z) = \frac{p(z)}{q(z)} = \frac{\Psi(\varsigma)}{q(\varsigma)s(\varsigma)}(z_1 A_1 + z_2 A_2 + \ldots + z_d A_d)\frac{\Psi(z)^T}{q(z)s(z)}.$$

Hence we get

$$\operatorname{Im} f(z) = \sum_{k=1}^{d} \operatorname{Im} z_k \Phi_k(z, \overline{z}), \qquad \operatorname{Re} f(z) = \sum_{k=1}^{d} \operatorname{Re} z_k \Phi_k(z, \overline{z}), \qquad (8.1)$$

$$W_{z_1}[q(z), p(z)] = \left(\frac{H(z)}{s(z)}\right)\left(\frac{H(z)}{s(z)}\right)^T, \qquad (8.2)$$

where $H(z) = \Psi(z) A_1^{1/2} = (h_1(z), \cdots, h_N(z))$ and

$$\Phi_1(z, \overline{z}) = \frac{1}{|q(z)|^2}\left(\frac{H(z)}{s(z)}\right)\left(\frac{H(z)}{s(z)}\right)^*; \quad \Phi_k(z, \overline{z}) = \frac{1}{|q(z)|^2}\left(\frac{\Psi(z)}{s(z)}\right) A_k \left(\frac{\Psi(z)}{s(z)}\right)^*, \quad k = 2, \cdots, d. \quad (8.3)$$

*The following observation is essential. In relation* (8.2), *the right side is a polynomial (we cancel* $s(z)^2$ *out of that fraction). If* $z \in \mathbb{R}^d$, *then the expressions*

$$\left(\frac{H(z)}{s(z)}\right)\left(\frac{H(z)}{s(z)}\right)^* \quad and \quad \left(\frac{\Psi(z)}{s(z)}\right) A_k \left(\frac{\Psi(z)}{s(z)}\right)^*$$

*in* (8.3) *are still polynomials, but if* $z \notin \mathbb{R}^d$, *then square of the modulus* $|s(z)|^2$ *no longer "cancels out". Then the complex zeros of* $s(z)$ *can be singularity points for the functions* $\operatorname{Re} f(z)$ *and* $\operatorname{Im} f(z)$ *in* (8.1).

Each irreducible factor $s_j(z)$ cannot be a divisor of all elements of the row vector $H(z)$. If $s_j(z)$ is a divisor of all elements $H(z)$, then it follows from (8.2) that $\hat{s}(z) = s(z)/s_j(z)$ is also Artin's denominator of the Wronskian $W_{z_1}[q(z), p(z)]$. This contradicts the minimality of $s(z)$. Let $\hat{h}(z)$ denote an element of $H(z)$ for which $s_j(z)$ is not a divisor.

Since $\overline{A}_k = A_k = A_k^T$, $k = 2, \cdots, d$, we see that there exist real matrices $R_k$ (generally rectangular) such that $A_k = R_k J_k R_k^*$, where $J_k$ are diagonal matrices with elements $\pm 1$ on the main diagonal. Let us introduce the designate $\Psi(z) R_k = (g_1^{(k)}(z), \cdots, g_{r_k}^{(k)}(z))$. Then

$$\Phi_1(z,\bar{z}) = \sum_{j=1}^{N} \frac{|h_j(z)|^2}{|q(z)s(z)|^2}; \quad \Phi_k(z,\bar{z}) = \sum_{j=1}^{n_k} \delta_j^{(k)} \frac{|g_j^{(k)}(z)|^2}{|q(z)s(z)|^2}, \quad \delta_j^{(k)} = \pm 1, \; k = 2,\cdots,d \; . \quad (8.4)$$

***Case* (a)**. Let $s_0(z)$, $\partial s_0(z)/\partial z_1 \neq 0$ not be a divisor of an element $\hat{h}(z)$ from $H(z)$. By Lemma 8.2, there exists a point $z' = (\eta_1', \hat{x}') \in \overline{\mathcal{H}}^d$, $\mathrm{Im}\,\eta_1' > 0$ such that $s(z') = s_0(z') = 0$, $\hat{h}(z') \neq 0$. It follows from (8.1) that

$$\mathrm{Im}\, f(z_1, \hat{x}') = \mathrm{Im}\, z_1 \cdot \Phi_1 = \mathrm{Im}\, z_1 \cdot \sum_{j=1}^{N} \frac{|h_j(z_1, \hat{x}')|^2}{|q(z_1, \hat{x}')s(z_1, \hat{x}')|^2}.$$

We obtain $\lim_{z_1 \to \eta_1'} \mathrm{Im}\, f(z_1, \hat{x}') = +\infty$. That is, *the function $f(z)$ has a singularity at the point* $(\eta_1', x_2', \cdots, x_d') \in \overline{\mathcal{H}}^d$, $\mathrm{Im}\,\eta_1' > 0$.

***Case* (b)**. Let Artin's denominator $s(z) = s(z_2, z_3, \cdots, z_d)$ be independent of the variable $z_1$. Without loss of generality, we can assume that $\partial s_0(z)/\partial z_2 \neq 0$, where $s_0(z_2, \cdots, z_d)$ is an irreducible factor of $s(z)$. Let $s_0(z)$ not be a divisor of an element $\hat{h}(z)$ from $H(z)$. By Lemma 8.2, there exists $\hat{z}' = (\eta_2', \hat{x}') \in \overline{\mathcal{H}}^{d-1}$, $\mathrm{Im}\,\eta_2' > 0$ such that $s_0(\hat{z}') = 0$. We have $\hat{h}(z_1', \eta_2', \hat{x}') \neq 0$ for almost all $z_1' \in \mathbb{C}$, $\mathrm{Im}\, z_1' > 0$ except for a finite number of points.

Without loss of generality, we can assume that $\mathrm{Re}\,\eta_2' \neq 0$. If $\mathrm{Re}\,\eta_2' = 0$, then instead of the variable $z_2$ consider a new variable $\hat{z}_2 = z_2 - x_2'$, $x_2' \in \mathbb{R}\setminus\{0\}$.

Let us show that for all $z_1' \in \mathbb{C}, \mathrm{Im}\, z_1' > 0$, the double limit

$$\lim_{\substack{z_1 \to z_1' \\ z_2 \to \eta_2'}} \Phi_1(z_1, z_2, \hat{x}') = \lim_{\substack{z_1 \to z_1' \\ z_2 \to \eta_2'}} \frac{1}{|q(z_1, z_2, \hat{x}')|^2} \sum_{j=1}^{N} \left|\frac{h_j(z_1, z_2, \hat{x}')}{s(z_2, \hat{x}')}\right|^2 \quad (8.5)$$

either does not exist or is equal $+\infty$.

For each $z_1' \in \mathbb{C}$, $\mathrm{Im}\, z_1' > 0$ there are two possibilities: either $\hat{h}(z_1', \eta_2', \hat{x}') = 0$, or $\hat{h}(z_1', \eta_2', \hat{x}') \neq 0$. In the first case, by Theorem 2.2, the limit (8.5) does not exist, and in the second case (8.5) is equal to $+\infty$.

Further reasoning is based on the relations

$$\mathrm{Im}\, f(z_1, z_2, \hat{x}') = \mathrm{Im}\, z_1 \cdot \Phi_1(z_1, z_2, \hat{x}') + \mathrm{Im}\, z_2 \cdot \Phi_2(z_1, z_2, \hat{x}'), \quad (8.6)$$

$$\mathrm{Re}\, f(z_1, z_2, \hat{x}') = \mathrm{Re}\, z_1 \cdot \Phi_1(z_1, z_2, \hat{x}') + \mathrm{Re}\, z_2 \cdot \Phi_2(z_1, z_2, \hat{x}') + \sum_{k=3}^{d} x_k' \Phi_k(z_1, z_2, \hat{x}'). \quad (8.7)$$

There are two alternative possibilities for the sum $\sum_{k=3}^{d} x_k' \Phi_k(z_1, z_2, \hat{x}')$ in (8.7):

**(b.1)**. For each $z_1' \in \mathbb{C}$, $\mathrm{Im}\, z_1' > 0$, $\mathrm{Re}\, z_1' \neq 0$ double limit

$$\lim_{\substack{z_1 \to z_1' \\ z_2 \to \eta_2'}} \sum_{k=3}^{d} x_k' \Phi_k(z_1, z_2, \hat{x}') \quad (8.8)$$

either does not exist or is equal to $\pm\infty$;

**(b.2).** There exists a finite limit (8.8) at the point $z_1' \in \mathbb{C}$, $\operatorname{Im} z_1' > 0$, $\operatorname{Re} z_1' \neq 0$.

*Case* **(b.1).** Limit (8.8) either does not exist or is equal to $\pm \infty$. For $\operatorname{Re} \eta_2' \neq 0$ there exists a point $(z_1', \eta_2', \hat{x}')$ such that $\operatorname{Im} z_1' > 0$, $\operatorname{Re} z_1' \neq 0$ and $\operatorname{Re} z_1' = \dfrac{\operatorname{Re} \eta_2'}{\operatorname{Im} \eta_2'} \operatorname{Im} z_1'$.

Suppose that the function $f(z)$ is regular at the point $(z_1', \eta_2', \hat{x}')$; then there exist finite limits

$$\lim_{\substack{z_1 \to z_1' \\ z_2 \to \eta_2'}} \operatorname{Im} f(z_1, z_2, \hat{x}') = \lim_{\substack{z_1 \to z_1' \\ z_2 \to \eta_2'}} \left[ \operatorname{Im} z_1 \cdot \Phi_1(z_1, z_2, \hat{x}') + \operatorname{Im} z_2 \cdot \Phi_2(z_1, z_2, \hat{x}') \right],$$

$$\lim_{\substack{z_1 \to z_1' \\ z_2 \to \eta_2'}} \operatorname{Re} f(z_1, z_2, \hat{x}') = \lim_{\substack{z_1 \to z_1' \\ z_2 \to \eta_2'}} \left[ \operatorname{Re} z_1 \cdot \Phi_1(z_1, z_2, \hat{x}') + \operatorname{Re} z_2 \cdot \Phi_2(z_1, z_2, \hat{x}') + \sum_{k=3}^d x_k' \Phi_k(z_1, z_2, \hat{x}') \right].$$

Since limit (8.8) either does not exist or is equal to $\pm \infty$, we see that

$$\lim_{\substack{z_1 \to z_1' \\ z_2 \to \eta_2'}} \left[ \operatorname{Re} z_1 \cdot \Phi_1(z_1, z_2, \hat{x}') + \operatorname{Re} z_2 \cdot \Phi_2(z_1, z_2, \hat{x}') \right]$$

must have the same property. From condition $\operatorname{Re} z_1' = K \cdot \operatorname{Im} z_1'$, where $K = (\operatorname{Re} \eta_2' / \operatorname{Im} \eta_2') \neq 0$ we obtain

$$\lim_{\substack{z_1 \to z_1' \\ z_2 \to \eta_2'}} \operatorname{Im} f(z_1, z_2, \hat{x}') = \lim_{\substack{z_1 \to z_1' \\ z_2 \to \eta_2'}} \left[ \operatorname{Im} z_1 \cdot \Phi_1(z_1, z_2, \hat{x}') + \operatorname{Im} z_2 \cdot \Phi_2(z_1, z_2, \hat{x}') \right] =$$

$$= (1/K) \cdot \lim_{\substack{z_1 \to z_1' \\ z_2 \to \eta_2'}} \left[ \operatorname{Re} z_1 \cdot \Phi_1(z_1, z_2, \hat{x}') + \operatorname{Re} z_2 \cdot \Phi_2(z_1, z_2, \hat{x}') \right].$$

The last limit either does not exist or is equal $\pm \infty$, which contradicts the assumption. *The function $f(z)$ has a singularity at the point $(z_1', \eta_2', \hat{x}') \in \overline{\mathcal{H}}^d$, $\operatorname{Im} z_1' > 0$.*

*Case* **(b.2).** Let us prove that $f(z_1, z_2, \hat{x}')$ also has a singularity at the point $(z_1', \eta_2', \hat{x}')$.

Assume the converse. Then function $f_{\hat{x}'}(z_1, z_2) = f(z_1, z_2, \hat{x}')$ is continuous at the point $(z_1', \eta_2')$, i.e., there exist finite limits

$$\lim_{\substack{z_1 \to z_1' \\ z_2 \to \eta_2'}} \operatorname{Re} f(z_1, z_2, \hat{x}'), \qquad \lim_{\substack{z_1 \to z_1' \\ z_2 \to \eta_2'}} \operatorname{Im} f(z_1, z_2, \hat{x}').$$

Then rational function $f_{\hat{x}'}(z_1, z_2)$ is continuous in a sufficiently small neighborhood of the point $(z_1', \eta_2')$. If $\lim_{\substack{z_1 \to z_1' \\ z_2 \to \eta_2'}} \sum_{k=3}^d x_k' \Phi_k(z_1, z_2, \hat{x}') = A \neq \pm \infty$, then (see (8.7), (8.6)) there exist finite limits

$$\lim_{\substack{z_1 \to z_1' \\ z_2 \to \eta_2'}} \left[ \operatorname{Re} z_1 \cdot \Phi_1(z_1, z_2, \hat{x}') + \operatorname{Re} z_2 \cdot \Phi_2(z_1, z_2, \hat{x}') \right], \tag{8.9}$$

$$\lim_{\substack{z_1 \to z_1' \\ z_2 \to \eta_2'}} \left[ \operatorname{Im} z_1 \cdot \Phi_1(z_1, z_2, \hat{x}') + \operatorname{Im} z_2 \cdot \Phi_2(z_1, z_2, \hat{x}') \right]. \tag{8.10}$$

Since the double limit (8.5) either does not exist or is equal to $+\infty$, we see that the limits of the first terms in (8.9) and (8.10) will be the same. By assumption, both limits (8.9) and (8.10) are finite. Then the points $z_1'$ and $\eta_2'$ must necessarily satisfy the relations

$$\operatorname{Re} z_1' = K \cdot \operatorname{Im} z_1', \quad \operatorname{Re} \eta_2' = K \cdot \operatorname{Im} \eta_2', \qquad (8.11)$$

where $K = \operatorname{Re} \eta_2' / \operatorname{Im} \eta_2' \neq 0$.

Condition (8.11) is violated at point $(z_1'', \eta_2')$, where $z_1'' = z_1' + \varepsilon$, $\varepsilon \in \mathbb{R} \setminus \{0\}$. For each such point $(z_1'', \eta_2')$, either limit (8.9) or (8.10) does not exist or is equal to $\pm\infty$. This contradicts the continuity of the function $f_{\hat{x}'}(z_1, z_2)$ in a small neighborhood of the point $(z_1', \eta_2')$. The function $f(z)$ has a singularity at the point $(z_1', \eta_2', \hat{x}') \in \overline{\mathcal{H}}^d$, $\operatorname{Im} z_1' > 0$. ∎

**8.3. Corollary.** *Proposition 8.1 remains valid if replace scalar form $p(z)$ by $m \times m$ matrix form $P(z) = P(z)^T$.*

In the proof, one should consider the diagonal elements $f_{ll}(z)$ of the matrix-valued function $f(z) = P(z)/q(z)$. ∎

**8.4. Sum of Squares Theorem.** *If $P(z)/q(z) \in \mathbb{RP}_d^{m \times m}$, then the partial Wronskians*

$$W_{z_k}[q, P] = q(z) \frac{\partial P(z)}{\partial z_k} - P(z) \frac{\partial q(z)}{\partial z_k}, \quad k = 1, \cdots, d \qquad (8.12)$$

*are SOS forms.*

**Proof.** By Proposition 3.1, Wronskians $W_{z_k}[q, P]$ are PSD forms. If $d \leq 2$, then each PSD form is a SOS form. Suppose $f(z) = P(z)/q(z) \in \mathbb{RP}_d^{m \times m}$, where $d \geq 3$. By Proposition 3.2, if $\deg_{z_k} P(z) > \deg_{z_k} q(z)$, then

$$f(z) = z_1 A_1 + \ldots + z_d A_d + f_1(z),$$

where $f_1(z) = P_1(z)/q(z) \in \mathbb{RP}_d^{m \times m}$, $A_k \geq 0$ and $\deg_{z_k} P_1(z) = \deg_{z_k} q(z)$, $k = 1, \cdots, d$.

If $W_{z_k}[q(z), P_1(z)]$ are SOS forms and $A_k \geq 0$, then Wronskians

$$W_{z_k}[q, P] = q(z)^2 A_k + W_{z_k}[q(z), P_1(z)]$$

will also be SOS forms.

Let now $\deg_{z_k} P(z) = \deg_{z_k} q(z)$, $k = 1, \cdots, d$.

Assume the converse. Then at least one of the Wronskians, say $W_{z_1}[q(z), P(z)]$, is a PSD not SOS form. By Proposition 8.1, there exists a point $z' = (z_1, z_2, \cdots, z_d) \in \overline{\mathcal{H}}^d$, $\operatorname{Im} z_1 > 0$ such that $f(z)$ has a singularity at this point. Since $f(z) \in \mathbb{RP}_d^{m \times m}$, we see that at each point $(z_1, z_2, \cdots, z_d) \in \overline{\mathcal{H}}^d$, $\operatorname{Im} z_1 > 0$ the inequality $\operatorname{Im} f(z_1, z_2, \cdots, z_d) \geq 0$. Then function $f(z)$ is regular at each such point. This contradiction proves the theorem. ∎

## 9. Long Resolvent Representations of Positive Real Functions

The following theorem and Corollary 9.2 are the main results of this article.

**9.1. Theorem.** *Let $f(z)$ be a rational function of class $\mathbb{RP}_d^{m \times m}$. Then:*
   (a). *There exists an $(m+l) \times (m+l)$ matrix pencil*

$$A(z) = z_1 A_1 + \ldots + z_d A_d = \begin{pmatrix} A_{11}(z) & A_{12}(z) \\ A_{21}(z) & A_{22}(z) \end{pmatrix} \qquad (9.1)$$

with are real positive semidefinite matrices $A_k \geq 0$, $k = 1, \cdots, d$;

(b). *The function $f(z)$ has the long-resolvent representation*

$$f(z) = A_{11}(z) - A_{12}(z) A_{22}(z)^{-1} A_{21}(z). \qquad (9.2)$$

**9.2. Corollary**. $\mathbb{R}\mathcal{P}_d^{m \times m} = \mathbb{R}\mathcal{B}_d^{m \times m}$ for each $d \geq 1$. □

The proof of Theorem 9.1 is based on the following generalization of Darlington's theorem to the case of a function of several variables. Recall that a rational function $P(z)/q(z)$ is called *multiaffine* if $\deg_{z_k} P(z) = \deg_{z_k} q(z) = 1$, $k = 1, \cdots, d$.

**9.3. Proposition**. Let $f(z) = \dfrac{P(z)}{q(z)} = \dfrac{z_1 P_1(\hat{z}) + P_2(\hat{z})}{z_1 q_1(\hat{z}) + q_2(\hat{z})} \in \mathbb{R}\mathcal{P}_d^{m \times m}$ be a multiaffine function, let $W_{z_1}[q, P] = \Phi_1(\hat{z}) \Phi_1^T(\hat{z})$ be a SOS form. If $\Phi_1(\hat{z})$, $\hat{z} = (z_2, \cdots, z_d)$ has the size $m \times r$, then

$$g(\hat{z}) = \begin{pmatrix} g_{11}(\hat{z}) & g_{12}(\hat{z}) \\ g_{12}(\hat{z})^T & g_{22}(\hat{z}) \end{pmatrix} = \begin{pmatrix} \dfrac{P_1(\hat{z})}{q_1(\hat{z})} & \dfrac{\Phi_1(\hat{z})}{q_1(\hat{z})} \\ \dfrac{\Phi_1^T(\hat{z})}{q_1(\hat{z})} & \dfrac{q_2(\hat{z})}{q_1(\hat{z})} I_r \end{pmatrix} \qquad (9.3)$$

*is a multiaffine function of class $\mathbb{R}\mathcal{P}_{d-1}^{(m+r) \times (m+r)}$, and*

$$f(z) = g_{11}(\hat{z}) - g_{12}(\hat{z}) \big( g_{22}(\hat{z}) + z_1 I_r \big)^{-1} g_{12}(\hat{z})^T. \qquad (9.4)$$

*Proof*. Representation (9.4) follows from the obvious identity

$$f(z) = \frac{z_1 P_1(\hat{z}) + P_2(\hat{z})}{z_1 q_1(\hat{z}) + q_2(\hat{z})} = \frac{P_1}{q_1} - \frac{\Phi_1 \Phi_1^T}{q_1^2 (z_1 + q_2/q_1)}. \qquad (9.5)$$

Let us prove $g(\hat{z}) \in \mathbb{R}\mathcal{P}_{d-1}^{(m+r) \times (m+r)}$. The multiaffinity of $g(\hat{z})$ is obvious. According to Theorem 3.7 (the criterion of positivity), it suffices to prove that the forms $F_k(\hat{z}) = q_1^2 \partial g(\hat{z})/\partial z_k$, $k = 2, \cdots, d$ are PSD forms. Since $f(z)$ is multiaffine, we have

$$f(z) = \frac{P(z)}{q(z)} = \frac{z_1 z_k \hat{P}_1 + z_1 \hat{P}_2 + z_k \hat{P}_3 + \hat{P}_4}{z_1 z_k \hat{q}_1 + z_1 \hat{q}_2 + z_k \hat{q}_3 + \hat{q}_4}. \qquad (9.6)$$

It follows from (9.6) that

$$F_k(\hat{z}) = q_1^2 \frac{\partial g(\hat{z})}{\partial z_k} = \begin{pmatrix} \hat{P}_1 \hat{q}_2 - \hat{P}_2 \hat{q}_1 & \Phi_k(\hat{z}) \\ \Phi_k^T(\hat{z}) & (\hat{q}_2 \hat{q}_3 - \hat{q}_1 \hat{q}_4) I_r \end{pmatrix},$$

where $\Phi_k(\hat{z}) = (z_k \hat{q}_1 + \hat{q}_2) \dfrac{\partial \Phi_1}{\partial z_k} - \hat{q}_1 \Phi_1$ does not depend on the variables $z_1$ and $z_k$.

Let us prove the identity

$$(\hat{P}_1 \hat{q}_2 - \hat{P}_2 \hat{q}_1)(\hat{q}_2 \hat{q}_3 - \hat{q}_1 \hat{q}_4) = \Phi_k(\hat{z}) \Phi_k^T(\hat{z}). \qquad (9.7)$$

We obtain

$$\Phi_k(\hat{z})\Phi_k^T(\hat{z}) = \left((z_k\hat{q}_1 + \hat{q}_2)\frac{\partial \Phi_1}{\partial z_k} - \hat{q}_1\Phi_1\right)\left((z_k\hat{q}_1 + \hat{q}_2)\frac{\partial \Phi_1^T}{\partial z_k} - \hat{q}_1\Phi_1^T\right) =$$

$$= (z_k\hat{q}_1 + \hat{q}_2)^2 \frac{\partial \Phi_1}{\partial z_k}\frac{\partial \Phi_1^T}{\partial z_k} - (z_k\hat{q}_1^2 + \hat{q}_1\hat{q}_2)\left(\frac{\partial \Phi_1}{\partial z_k}\Phi_1^T + \Phi_1\frac{\partial \Phi_1^T}{\partial z_k}\right) + \hat{q}_1^2\Phi_1\Phi_1^T. \quad (9.8)$$

It follows from (9.6) that

$$\Phi_1(\hat{z})\Phi_1^T(\hat{z}) = q(z)^2 \frac{\partial f(z)}{\partial z_1} = z_k^2(\hat{q}_3\hat{P}_1 - \hat{q}_1\hat{P}_3) + z_k(\hat{q}_4\hat{P}_1 - \hat{q}_1\hat{P}_4 + \hat{q}_3\hat{P}_2 - \hat{q}_2\hat{P}_3) + (\hat{q}_4\hat{P}_2 - \hat{q}_2\hat{P}_4).$$

Hence

$$\left(\frac{\partial \Phi_1}{\partial z_k}\Phi_1^T + \Phi_1\frac{\partial \Phi_1^T}{\partial z_k}\right) = 2z_k(\hat{q}_3\hat{P}_1 - \hat{q}_1\hat{P}_3) + (\hat{q}_4\hat{P}_1 - \hat{q}_1\hat{P}_4 + \hat{q}_3\hat{P}_2 - \hat{q}_2\hat{P}_3),$$

$$\frac{\partial \Phi_1}{\partial z_k}\frac{\partial \Phi_1^T}{\partial z_k} = (\hat{q}_3\hat{P}_1 - \hat{q}_1\hat{P}_3).$$

Substituting the last 3 relations in (9.8), we see that identity (9.7) holds.

Since $h = \dfrac{q}{\partial q / \partial z_1}\bigg|_{z_1=0} = \left(z_1 + \dfrac{z_k\hat{q}_3 + \hat{q}_4}{z_k\hat{q}_1 + \hat{q}_2}\right)_{z_1=0} = \dfrac{z_k\hat{q}_3 + \hat{q}_4}{z_k\hat{q}_1 + \hat{q}_2} \in \mathbb{R}\mathcal{P}_{d-1}^{1\times 1}$, we see that $(\hat{q}_2\hat{q}_3 - \hat{q}_1\hat{q}_4)$ is PSD form. Then at $\hat{x} \in \mathbb{R}^{d-2}$ we obtain

$$F_k(\hat{x}) = \begin{pmatrix} I_m & \Phi_k(\hat{q}_2\hat{q}_3 - \hat{q}_1\hat{q}_4)^{-1} \\ 0 & I_r \end{pmatrix}\begin{pmatrix} 0 & 0 \\ 0 & (\hat{q}_2\hat{q}_3 - \hat{q}_1\hat{q}_4)I_r \end{pmatrix}\begin{pmatrix} I_m & 0 \\ (\hat{q}_2\hat{q}_3 - \hat{q}_1\hat{q}_4)^{-1}\Phi_k^T & I_r \end{pmatrix} \geq 0. \blacksquare$$

**9.4. Lemma**. *Suppose* $f(z) = g_{11} - g_{12}(g_{22} + z_1 I_{r_1})^{-1}g_{12}^T$, *where*

$$\begin{pmatrix} g_{11} & g_{12} \\ g_{12}^T & g_{22} \end{pmatrix} = \begin{pmatrix} a_{11} & a_{12} \\ a_{12}^T & a_{22} \end{pmatrix} - \begin{pmatrix} a_{13} \\ a_{23} \end{pmatrix}\left[(d + z_2)I_{r_2}\right]^{-1}\begin{pmatrix} a_{13}^T & a_{23}^T \end{pmatrix}.$$

*Then*

$$f(z) = A_{11}(z) - \begin{pmatrix} A_{12}(z) & A_{13}(z) \end{pmatrix}\begin{pmatrix} A_{22}(z) & A_{23}(z) \\ A_{23}^T(z) & A_{33}(z) \end{pmatrix}^{-1}\begin{pmatrix} A_{12}^T(z) \\ A_{13}^T(z) \end{pmatrix}, \quad (9.9)$$

*where*

$$\begin{pmatrix} A_{11}(z) & A_{12}(z) & A_{13}(z) \\ A_{12}^T(z) & A_{22}(z) & A_{23}(z) \\ A_{13}^T(z) & A_{23}^T(z) & A_{33}(z) \end{pmatrix} = \begin{pmatrix} a_{11} & a_{12} & a_{13} \\ a_{12}^T & a_{22} & a_{23} \\ a_{13}^T & a_{23}^T & d \cdot I_{r_2} \end{pmatrix} + z_1\begin{pmatrix} 0 & 0 & 0 \\ 0 & I_{r_1} & 0 \\ 0 & 0 & 0 \end{pmatrix} + z_2\begin{pmatrix} 0 & 0 & 0 \\ 0 & 0 & 0 \\ 0 & 0 & I_{r_2} \end{pmatrix}.$$

***Proof***. The identity (9.9) is verified by direct computation. $\blacksquare$

***Proof of Theorem 9.1***. Suppose $f(z) = P(z)/q(z) \in \mathbb{R}\mathcal{P}_d^{m\times m}$, $\deg_{z_k} f(z) = n_k$, $k = 1, \cdots, d$. Without loss of generality (Proposition 3.2), we can consider $\deg_{z_k} P(z) = \deg_{z_k} q(z)$, $k = 1, \cdots, d$. Applying to the function $f(z)$ Degree Reduction Operator $\mathrm{D}_{z_k}^{n_k}$ in each variable $z_k$, $k = 1, \cdots, d$, we obtain multiaffine function

$$\hat{f}(\varsigma_1,\cdots,\varsigma_n) = \hat{P}(\varsigma_1,\cdots,\varsigma_n)/\hat{q}(\varsigma_1,\cdots,\varsigma_n) \qquad (9.10)$$

in variables $\varsigma_1,\cdots,\varsigma_n$, $n = n_1 + \ldots + n_d$. It follows from Theorem 3.4 that $\hat{f}(\varsigma_1,\cdots,\varsigma_n)$ is a positive real function. Moreover, $\deg_{\varsigma_j} \hat{P}(\varsigma_1,\cdots,\varsigma_n) = \deg_{\varsigma_j} \hat{q}(\varsigma_1,\cdots,\varsigma_n) = 1$, $j = 1,\cdots,n$. Then

$$\hat{f}(\varsigma) = \frac{P(\varsigma_1,\cdots,\varsigma_n)}{q(\varsigma_1,\cdots,\varsigma_n)} = \frac{\varsigma_1 P_1(\hat{\varsigma}) + P_2(\hat{\varsigma})}{\varsigma_1 q_1(\hat{\varsigma}) + q_2(\hat{\varsigma})},$$

where $q_1(\hat{\varsigma}) \not\equiv 0$. By Theorem 8.4, there exists a multiaffine $m \times r_1$ matrix-valued form $\Phi_1(\hat{\varsigma})$ such that

$$W_{\varsigma_1}[q,P] = P_1(\hat{\varsigma})q_2(\hat{\varsigma}) - P_2(\hat{\varsigma})q_1(\hat{\varsigma}) = \Phi_1(\hat{\varsigma})\Phi_1^T(\hat{\varsigma}).$$

It follows from Proposition 9.3 that function

$$g^{(1)}(\hat{\varsigma}) = \begin{pmatrix} g_{11}^{(1)}(\hat{\varsigma}) & g_{12}^{(1)}(\hat{\varsigma}) \\ g_{12}^{(1)}(\hat{\varsigma})^T & g_{22}^{(1)}(\hat{\varsigma}) \end{pmatrix} = \begin{pmatrix} \dfrac{P_1(\hat{\varsigma})}{q_1(\hat{\varsigma})} & \dfrac{\Phi_1(\hat{\varsigma})}{q_1(\hat{z})} \\ \dfrac{\Phi_1^T(\hat{\varsigma})}{q_1(\hat{z})} & \dfrac{q_2(\hat{\varsigma})}{q_1(\hat{z})} I_{r_1} \end{pmatrix}$$

belongs to the class $\mathbb{RP}_{n-1}^{(m \times r_1) \times (m \times r_1)}$, and

$$\hat{f}(\varsigma_1,\cdots,\varsigma_n) = g_{11}^{(1)}(\hat{\varsigma}) - g_{12}^{(1)}(\hat{\varsigma})\left(g_{22}^{(1)}(\hat{\varsigma}) + \varsigma_1 I_{r_1}\right)^{-1} g_{12}^{(1)}(\hat{\varsigma})^T.$$

The matrix-valued function $g^{(1)}(\hat{\varsigma})$ depends only on $(n-1)$ variables $\varsigma_2,\cdots,\varsigma_n$ and satisfies the conditions of Proposition 9.3. Then

$$\begin{pmatrix} g_{11}^{(1)}(\hat{\varsigma}) & g_{12}^{(1)}(\hat{\varsigma}) \\ g_{12}^{(1)}(\hat{\varsigma})^T & g_{22}^{(1)}(\hat{\varsigma}) \end{pmatrix} = \begin{pmatrix} g_{11}^{(2)} & g_{12}^{(2)} \\ g_{12}^{(2)T} & g_{22}^{(2)} \end{pmatrix} - \begin{pmatrix} g_{13}^{(2)} \\ g_{23}^{(2)} \end{pmatrix}\left[(g_{33}^{(2)} + \varsigma_2)I_{r_2}\right]^{-1}\begin{pmatrix} g_{13}^{(2)T} & g_{23}^{(2)T} \end{pmatrix}.$$

From Lemma 9.4 we obtain

$$\hat{f}(\varsigma_1,\cdots,\varsigma_n) = g_{11}^{(2)} - \begin{pmatrix} g_{12}^{(2)} & g_{13}^{(2)} \end{pmatrix}\left[\begin{pmatrix} g_{22}^{(2)} & g_{23}^{(2)} \\ g_{23}^{(2)T} & g_{33}^{(2)} I_{r_2} \end{pmatrix} + \begin{pmatrix} \varsigma_1 I_{r_1} & 0 \\ 0 & \varsigma_2 I_{r_2} \end{pmatrix}\right]^{-1}\begin{pmatrix} g_{12}^{(2)T} \\ g_{13}^{(2)T} \end{pmatrix},$$

where function

$$g^{(2)} = \begin{pmatrix} g_{11}^{(2)} & g_{12}^{(2)} & g_{13}^{(2)} \\ g_{12}^{(2)T} & g_{22}^{(2)} I_{r_1} & g_{23}^{(2)} \\ g_{13}^{(2)T} & g_{23}^{(2)T} & g_{33}^{(2)} I_{r_2} \end{pmatrix}$$

depends only on $(n-2)$ variables $\varsigma_3,\cdots,\varsigma_n$ and satisfies the conditions of Proposition 9.3.

Continuing the process, at the $(n-1)$ step we get a positive real matrix-valued function of one variable $\varsigma_n$:

$$g^{(n-1)}(\varsigma_n) = \begin{pmatrix} g_{11}^{(n-1)} & g_{12}^{(n-1)} & \cdots & g_{1n}^{(n-1)} \\ g_{12}^{(n-1)T} & g_{22}^{(n-1)} I_{r_1} & \cdots & g_{2n}^{(n-1)} \\ \vdots & \vdots & \ddots & \vdots \\ g_{1n}^{(n-1)T} & g_{2n}^{(n-1)T} & \cdots & g_{nn}^{(n-1)} I_{r_{n-1}} \end{pmatrix} = \varsigma_n \begin{pmatrix} A_{11} & A_{12} & \cdots & A_{1n} \\ A_{12}^T & A_{22} & \cdots & A_{2n} \\ \vdots & \vdots & \ddots & \vdots \\ A_{1n}^T & A_{2n}^T & \cdots & A_{nn} \end{pmatrix} = \varsigma_n A_n, \qquad (9.11)$$

where $A_n \geq 0$ is real positive semidefinite matrix. We have long-resolvent representation

$$\hat{f}(\varsigma_1,\cdots,\varsigma_n) = A_{11}(\varsigma) - A_{12}(\varsigma)A_{22}(\varsigma)^{-1}A_{11}(\varsigma)^T.$$

A positive real matrix pencil has the form

$$\begin{pmatrix} A_{11}(\varsigma) & A_{12}(\varsigma) \\ A_{12}(\varsigma)^T & A_{22}(\varsigma) \end{pmatrix} = \begin{pmatrix} 0 & 0 & \cdots & 0 \\ 0 & \varsigma_1 I_{r_1} & \cdots & 0 \\ \vdots & \vdots & \ddots & \vdots \\ 0 & 0 & \cdots & \varsigma_{n-1} I_{r_{n-1}} \end{pmatrix} + \varsigma_n \begin{pmatrix} A_{11} & A_{12} & \cdots & A_{1n} \\ A_{12}^T & A_{22} & \cdots & A_{2n} \\ \vdots & \vdots & \ddots & \vdots \\ A_{1n}^T & A_{2n}^T & \cdots & A_{nn} \end{pmatrix}.$$

The degree reduction operator (3.10) is invertible. If the variables of each $k$-th group are replaced by the variable $z_k$, then for the original positive function $f(z_1,\cdots,z_d)$ we obtain a long-resolvent representation with positive semidefinite matrices of pencil. ∎

## References.